\documentclass{amsart}

\usepackage[utf8]{inputenc} 

\usepackage{amsmath}
\usepackage{amsfonts}
\usepackage{amssymb}
\usepackage{amsthm}
\usepackage{commath}
\usepackage{mathrsfs}
\usepackage{mathtools}
\usepackage{bm}
\usepackage{enumitem}

\newtheorem{definition}{Definition}
\newtheorem{proposition}{Proposition}
\theoremstyle{plain}

\usepackage{graphicx}

\usepackage{xspace}
\usepackage{wrapfig}
\usepackage{subcaption}
\usepackage{float}
\floatstyle{plaintop}
\restylefloat{table}

\DeclareMathAlphabet{\mathpzc}{OT1}{pzc}{m}{it}

\usepackage[linesnumbered]{algorithm2e}
\usepackage[noend]{algpseudocode}
\usepackage[skins]{tcolorbox}
\SetKwInOut{Parameter}{Parameters}
\SetKwInOut{KwIn}{Input}
\SetKwInOut{KwOut}{Output}
\usepackage{varioref}

\usepackage{booktabs}
\usepackage{array}
\usepackage{hyperref}

\begin{document}

\title[A greedy algorithm for optimal heating]{A greedy algorithm for optimal heating in powder-bed-based additive manufacturing}

\author{Robert Forslund}
\address{Robert Forslund\\
Department of Mathematical Sciences\\
Chalmers University of Technology and University of Gothenburg\\
SE-412 96 Gothenburg\\
Sweden}
\email{robfor@chalmers.se}

\author{Anders Snis}
\address{Anders Snis\\
Arcam EBM\\
Krokslätts fabriker 27A\\
SE-431 37 Mölndal\\
Sweden}

\author[S.~Larsson]{Stig Larsson}
\address{Stig Larsson\\
Department of Mathematical Sciences\\
Chalmers University of Technology and University of Gothenburg\\
SE-412 96 Gothenburg\\
Sweden}


\begin{abstract} 
Powder-bed-based additive manufacturing involves melting of a powder bed using a moving laser or electron beam as a heat source. In this paper, we formulate an optimization scheme that aims to control this type of melting. The goal consists of tracking maximum temperatures on lines that run along the beam path. Time-dependent beam parameters (more specifically, beam power, spot size, and speed) act as control functions. The scheme is greedy in the sense that it exploits local properties of the melt pool in order to divide a large optimization problem into several small ones. As illustrated by numerical examples, the scheme can resolve heat conduction issues such as concentrated heat accumulation at turning points and non-uniform melt depths. 
\end{abstract}

\maketitle

\section{Introduction}
\label{sec:intro}
\noindent Powder bed fusion (PBF) is a type of additive manufacturing (AM) where metal powder is melted by a laser or electron beam in a layer-wise fashion to enable the production of geometrically complex parts \cite{heinl}. AM undergoes continuous progress towards a technology that is robust and efficient, but there are still issues when it comes to quality and repeatability.

It is important for completed parts to meet the mechanical requirements and quality standards specified by the applications for which they are manufactured. The qualities of a completed part, such as tensile strength and surface roughness, depend on the melting process, which in turn is governed by several dozen material and process parameters such as preheating temperature, powder packing ratio, and beam speed, among other. Correlations between process parameters, process signatures (such as melt pool size and temperature), and product qualities are presented in \cite{mani} and references therein. Due to these correlations, the design of process settings requires critical attention. However, this involves extensive and costly experimental work and the resulting melting schemes implemented in machines rely on an excessive amount of parameters and functions in order to account for the dynamics of the melting process. This approach makes it difficult to optimize PBF and limits both the number of applications and the number of materials available for manufacturing. 

As of late, computation based techniques are used to improve PBF related process control, i.e., to tackle the question of how to select different process parameters in order to build parts with desired properties \cite{mani, kingetal}. A common approach is design of (computational) experiments (DoE) \cite{marklandkorner, zeng}, which is an exploratory tool used to identify parameters that influence the qualities of the completed part. In \cite{ma}, a DoE with a finite element model determines that the beam power and beam speed are the two process controllable parameters that have the largest impact on peak temperature during a single track melt. In \cite{kamath}, a DoE based on a simple thermal model \cite{eagerandtsai} is used to determine optimal process parameters for building high density parts. Also this study shows that the beam power and beam speed have the largest impact on melt pool width and depth. Empirical modelling techniques such as artificial neural networks are also used to determine the optimal selection of process parameters, see, e.g., \cite{ning} as well as \cite{garg} and references therein.

While DoE is useful for developing a deeper understanding of the melting process, it can be a tedious affair to use them for optimization. An optimal control approach might be better suited for that purpose, as it starts off with preset, desired melt characteristics and seeks corresponding optimal process parameters via a mathematical optimization problem. 

Attempts of optimizing melt pool characteristics have been made based on the two-phase Stefan problem, where the free boundary between solid and liquid is understood via the Stefan condition. In the context of PBF, the free boundary characterizes the size and shape of the melt pool. In \cite{HINZE2007657}, the two-phase Stefan problem in a container is considered, and the temperature on the container boundary is optimized with respect to a desired transient evolution of the free boundary. Optimization problems based on quasi-steady state formulations of the two-phase Stefan problem are solved in \cite{Volkov2009, cao}. Here the desired free boundary between liquid and melt is prescribed and the goal involves tracking the melt temperature on the prescribed free boundary.

This paper is part of an effort that aims to reduce the number of parameters needed in the design of melt schemes. We present a simulation based framework that seeks to facilitate process optimization and material development, while keeping computational costs at a minimum. Thus, computational efficiency is prioritized as we trade a certain level of detail for a fast optimization scheme that can be applied to the melting of large domains. The scheme is efficient for three reasons:
    \begin{itemize}
    \item[$\triangleright$] The continuum thermal model for describing heat conduction includes an analytic solution that allows for fast and pointwise computation of temperatures during melting \cite{john, forslund}. The model assumes that the beam parameters are piecewise constant in time. We remark that beam parameters in actual AM machines are often set to vary in such a discontinuous way. The model does not include phase change, and does in particular not capture the solid-liquid interface. Instead, we use the term melt pool to simply denote the region where the temperature is larger than the melt temperature.
    \item[$\triangleright$] The formulation of the optimization problem involves a severe model order reduction. Rather than striving for some desired temperature distribution that is difficult to express, the goal consists of tracking preset reference maximum temperatures on lines that run along the beam path.
    \item[$\triangleright$] The resulting optimization problem is solved by a greedy algorithm that divides it into several small problems that are easier to solve. These sub-problems are solved consecutively as the beam traverses the powder bed.
    \end{itemize}

\noindent Together, this framework comprises an optimization scheme for process control as suggested in \cite{snis2}.

As noted earlier, DoE suggests that the beam parameters have the largest impact on the temperature distribution during melting and, ultimately, on the quality of the completed part. For this reason, we choose the beam power, spot size and speed as control variables. Our scheme allows for time-dependent beam parameters, which increases their ability to control the melting process. An ability to optimize these beam parameters is useful not only for validation, but also in order to speed up the development of process settings for new (and old) materials. Since optimization is an iterative process, the solver of the forward problem needs to be highly efficient. The analytic solution provides such efficiency.

The remainder of this paper is organized as follows. Section \ref{sec:thermal_model} describes the thermal model and its corresponding analytic solution. The optimization problem is formulated in Section \ref{sec:optimization_problem}. Section \ref{sec:greedy_algorithm} and Section \ref{sec:alt_greedy_algorithm} propose two greedy algorithms for solving said problem. In Section \ref{sec:examples} we apply our optimization scheme in numerical examples. Here we also detail how the combination of the two greedy algorithms can aid in the development of so called beam parameter functions. Additional comments and concluding remarks are given in Section \ref{sec:conclusions}.


\section{Thermal model}
\label{sec:thermal_model}
\noindent Consider the heat equation on the lower half space $\Omega = \mathbb{R}^2 \times \mathbb{R}_-$ during a time span $ \mathcal{T} = (0, T]$. Let $\Gamma$ denote the surface boundary $z = 0$ and let $u_\mathrm{init}$ denote the constant initial temperature. The beam travels on the surface $\Gamma$ along a preset, piecewise linear path
\[\mathcal{C}^\mathrm{s} = \{\bm{x}^\mathrm{s}(t): t \in \mathcal{T}\},\]
\noindent where $\bm{x}^\mathrm{s}(t) = (x^\mathrm{s}(t), y^\mathrm{s}(t), 0)$ is the position of the center of the beam at time $t$. The heat flux $\Phi$ due to a scanning electron beam is modeled as a Gaussian function
\[\Phi = \Phi(x, y, t)  = \frac{P(t)}{2\pi\sigma(t)^2}\,\mathrm{exp} \Big(-\tfrac{\left(x-x^\mathrm{s}\left(t\right)\right)^2+\left(y-y^\mathrm{s}\left(t\right)\right)^2}{2\sigma(t)^2}\Big).\]
\noindent The three beam parameters are the absorbed beam power $P(t)$, the beam spot size $\sigma(t)$, and the beam speed $v(t) = |\bm{v}(t)|$. Here $\bm{v}(t) = (v_x(t), v_y(t), 0) = (v(t)\cos{\theta(t)}, v(t)\sin{\theta(t)}, 0)$, where $\theta(t)$ is the angle between the positive $x$-axis and the direction of the path. This angle is known for any $t$ since the beam path $\mathcal{C}^\mathrm{s}$ is preset and it follows that $v$ uniquely defines the vector $(v_x, v_y, 0)$. The position of the beam $\bm{x}^\mathrm{s}(t)$ depends on the speed with which the beam has traveled the path $\mathcal{C}^\mathrm{s}$ up to time $t$.  The beam parameters are often set to vary in a piecewise constant fashion in AM machines. The following definition makes the concept of piecewise constant beam parameters more precise.
\begin{definition}\label{def:pw}
Given times $0 = t_0 < t_1 < \hdots < t_N = T$, 
\[(t_{n-1}, t_n] = \big(t_n^\mathrm{i}, t_n^\mathrm{f}\big], \,\,\,\, n = 1, 2, \hdots, N,\]
is a partition of $\mathcal{T}$ consisting of $N$ segments. Index $n$ indicates the $n^{th}$ segment in the partition, and a segment in turn is a collection of the following data:
\begin{itemize}
\item[--] $t_n^\mathrm{i}, \,\,t_n^\mathrm{f}$; an initial time and final time, respectively,
\item[--] $(P_n, \sigma_n, v_n)$; a triplet of power, spot size, and speed such that 
\[\left(P(t), \sigma(t), v(t)\right) = (P_n, \sigma_n, v_n)    
\,\,\,\mathrm{if}\,\, t \in \big(t_n^\mathrm{i}, t_n^\mathrm{f}\big],
\]
\item[--] $\ell_n^\mathrm{s}$; a line traversed by the beam between times $t_n^\mathrm{i}$ and $t_n^\mathrm{f}$, 
\item[--] $\bm{x}_n^\mathrm{i}, \,\,\bm{x}_n^\mathrm{f}$; the initial position and final position of $\ell_n^\mathrm{s}$, respectively.
\item[--] $\theta_n  =\mathrm{tan}^{-1}\left(\frac{y_n^f - y_n^i}{x_n^f - x_n^i}\right)$; the angle between the positive $x$-axis and the direction of the path.
\end{itemize}
\end{definition}
\noindent With $\hat{z}$ the outward unit normal of $\Omega$, the heat transfer problem can be written as
\begin{equation}
\begin{aligned}
\rho c_p \frac{\partial u}{\partial t} - \nabla \cdot ( \lambda \nabla u) &= 0 & \mathrm{in}&  \,\,\,\,\Omega \times \mathcal{T},\\
(\lambda  \nabla u) \cdot \hat{z} &= \Phi & \mathrm{on}& \,\,\,\, \Gamma \times \mathcal{T},\\
u(\cdot, 0) &= u_\mathrm{init} & \mathrm{in} & \,\,\,\, \Omega,
\end{aligned}
\label{eq:pde}
\end{equation}
\noindent where $\rho$, $c_p$, $\lambda $ denote density, heat capacity, and thermal conductivity, respectively. These material parameters are assumed to be constant. This gives us
the thermal diffusivity $\kappa = \lambda/\rho c_p$. Problem \eqref{eq:pde} has an analytic solution \cite{john, forslund}. We refer to these sources for a detailed derivation of this solution and merely outline it here. 
\begin{proposition}\label{prop:analytic_solution}
Given a partition as in Definition \ref{def:pw}, the solution of problem \eqref{eq:pde} can be written as
\begin{align*}
u(\bm{x}, t) &= u_\mathrm{init} + \sum_{k=1}^{n-1} u_{n,k}^\mathrm{I}(\bm{x}, t) + u_n^\Phi(\bm{x}, t) \,\,\, \,\,\,\mathrm{for}\,\, t \in \big(t_n^\mathrm{i}, t_n^\mathrm{f}\big], \,\,\, n = 1, 2, \hdots, N,
\end{align*}
\noindent where $u_{n,k}^\mathrm{I}$ is the temperature due to the earlier scanning of segment $k<n$ and $u_n^\Phi$ is the temperature due to the current scanning of segment $n$. 
\end{proposition}

\noindent The analytic expressions of $u_{n,k}^\mathrm{I}$ and $u_{n}^\Phi$ are derived in \cite{forslund}, wherein it also described how the solution can be efficiently computed.

The power is largely determined by the beam current, and this current can not be adjusted at a fast rate. Therefore, $P$ is set to be constant for all $t \in \mathcal{T}$ for the remainder of this paper. It should be noted, however, that if one would be interested in optimizing all three beam parameters, the following extends to the case of non-constant power as well.

A remark on our thermal model is in order, as it does not include cooling, the latent heat of fusion nor a description of the solid-liquid intersection between powder and melt pool, which makes the notion of a melt pool quite fuzzy. Here we use the term melt pool to simply denote the volume where the temperature is larger than the melt temperature. Hence we use the isothermal $\{\bm{x}(t): u(\bm{x}, t) = u_\mathrm{melt}\}$, where $u_\mathrm{melt}$ is the melt temperature of the powder, to represent the solid-liquid interface. Furthermore, since the model is a continuum model, it breaks down on the mesoscale where we see phenomena such as balling, inter-capillary effects, Plateau--Rayleigh instabilities, thermal expansion, among other \cite{khairallah, kingkhairallah}. Despite these restrictions, it is anticipated that effective parameters, tuned via comparisons with experiments, can be used to make the thermal model reliable enough for control and optimization. The optimization problem described below also aligns with the overarching aim to reduce the number of parameters needed for process control. 


\section{Formulation of the optimization problem}
\label{sec:optimization_problem}
\noindent The goal is to optimize the melting process with respect to the beam spot size and beam speed. Since these beam parameters are defined in a piecewise constant fashion according to Definition \ref{def:pw}, we can write
\begin{align*}
\sigma(t) &= \sum_{k=1}^N \sigma_k \chi_{(t_k^\mathrm{i}, t_k^\mathrm{f}]},\\
v(t) &= \sum_{k=1}^N v_k \chi_{(t_k^\mathrm{i}, t_k^\mathrm{f}]},
\end{align*}
\noindent where $\chi$ is the indicator function. Since we aim to optimize the speed, either the times $t_n^\mathrm{i}$, $t_n^\mathrm{f}$ or the positions $\bm{x}_n^\mathrm{i}$, $\bm{x}_n^\mathrm{f}$, $n = 1, 2, \hdots, N$, in Definition \ref{def:pw} will have to vary during optimization. Due to the beam path being preset and reasons that will become clear in the next section, it is better to fix the positions. Therefore it is more appropriate to express $\sigma$ and $v$ as space dependent functions instead. To this end, introduce the scanning distance
\[
\gamma(\bm{x}^\mathrm{s}) = \sum_{k=1}^{n-1} |\bm{x}^\mathrm{f}_k - \bm{x}^\mathrm{i}_k| +  |\bm{x}^\mathrm{s} - \bm{x}^\mathrm{i}_n|, \,\, \mathrm{if} \,\, \bm{x}^\mathrm{s} \in \ell_n^s, \,\, n = 1, 2, \hdots, N.
\] 
\noindent Then we have the following beam parameter functions:
\begin{equation}
\begin{split}
\sigma(\gamma(\bm{x}^\mathrm{s})) &= \sum_{k=1}^N \sigma_k \chi_{(\gamma(\bm{x}_k^\mathrm{i}), \gamma(\bm{x}_k^\mathrm{f})]},\\
v(\gamma(\bm{x}^\mathrm{s})) &= \sum_{k=1}^N v_k \chi_{(\gamma(\bm{x}_k^\mathrm{i}), \gamma(\bm{x}_k^\mathrm{f})]}.
\end{split}
\label{eq:parameter_functions}
\end{equation}
\noindent From \eqref{eq:parameter_functions}, a decision vector can immediately be extracted as $\bm{d} =(\bm{\sigma}, \bm{v}) =  \{(\sigma_k, v_k)\}_{k=1}^N = (\sigma_1, v_1, \sigma_2, v_2, \hdots, \sigma_N, v_N)$. The variables in the decision vector are bounded due to practical limitations, $\bm{d}_\mathrm{min} \leq \bm{d} \leq \bm{d}_\mathrm{max}$.

The control of the melting process is a multiobjective optimization problem due to the many correlations between process parameters, process signatures, and product qualities. The qualities of the final part are strongly dependent on the temperatures obtained during the melting process \cite{mukherjee}. The characteristics of the melt pool are important process signatures. If the melt pool is too small relative to the line offset (i.e., the distance between two hatch lines) and layer depth, powder might be left unmelted between hatch lines or between layers, causing discontinuities and porosity. Furthermore, high surface temperatures might result in too much evaporation and subsequent recoil pressure, which can result in undesired material transport such as ejection of molten materials that later cause defects \cite{gong} or formations of small ridges that prohibit the deposition of new powder layers and thus cause the manufacturing process to terminate \cite{vandenbroucke}.  Qualitatively, therefore, the choice of cost functional can be motivated by the desire to 
\begin{enumerate}
    \item maintain a uniform and appropriately sized melt pool during melting, and
    \item avoid too high surface temperatures.
\end{enumerate}


\subsection{A reductive approach}
\label{sec:reduction}
\begin{figure}[h!]
\captionsetup{width=0.9\linewidth}
\includegraphics[width=0.9\linewidth]{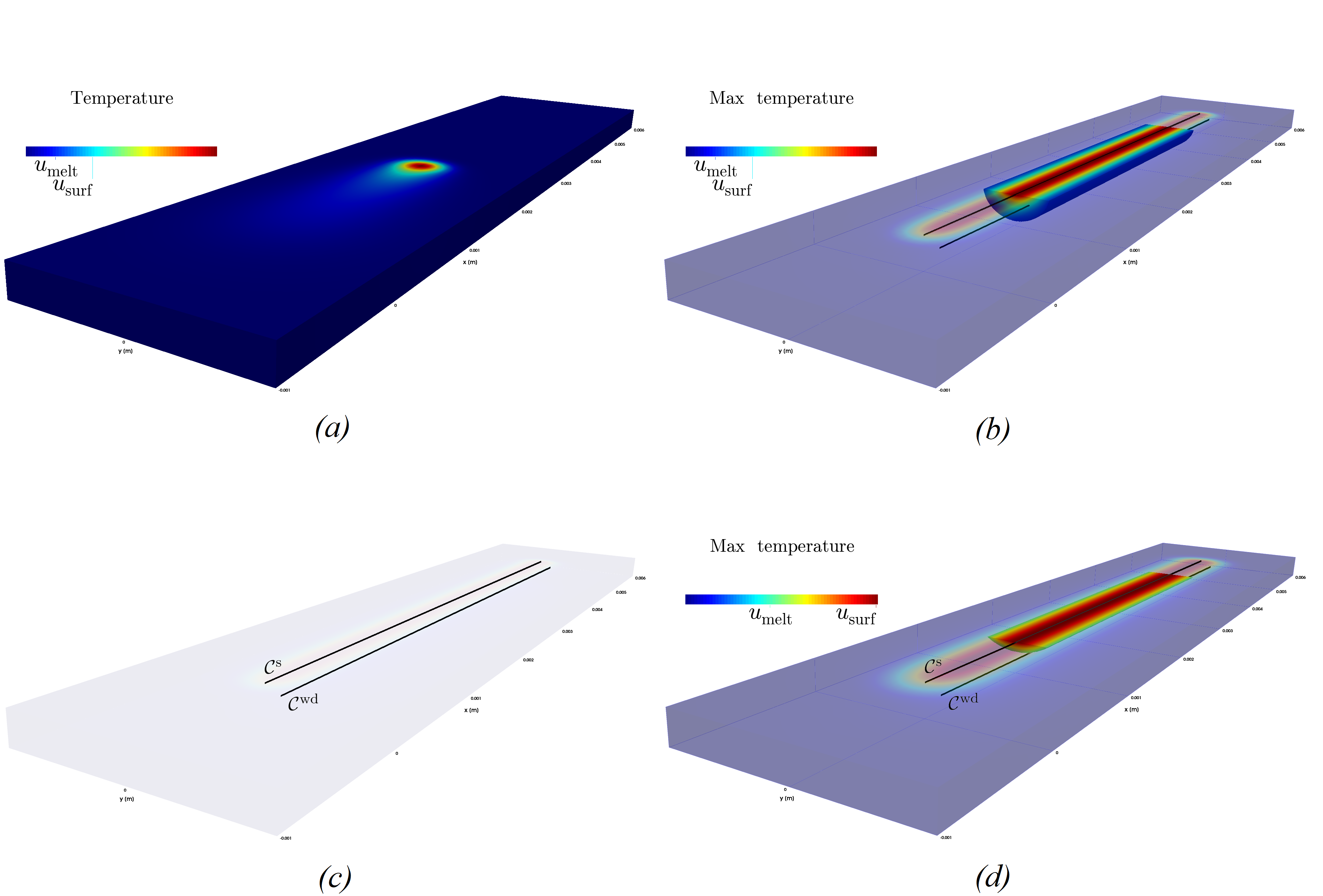}
\caption{The desire to optimize the size and shape of the melt pool is reduced to a problem of tracking maximum temperatures on paths. (a) Temperature distribution due to a moving beam. (b) Maximum temperature obtained during melting, with the solidified volume highlighted. (c) Introduction of beam path $\mathcal{C}^\mathrm{s}$ and secondary path $\mathcal{C}^\mathrm{wd}$. The secondary path is drawn along the desired liquid-solid interface. (d) Maximum temperature after optimization, with the solidified volume highlighted. The shape of this volume is optimized due to the tracking of reference temperatures $u_\mathrm{surf}$ and $u_\mathrm{melt}$ on $\mathcal{C}^\mathrm{s}$ and $\mathcal{C}^\mathrm{wd}$, respectively.}
\label{fig:reduction}
\end{figure}

\noindent Consider a beam scanning along a path $\mathcal{C}^\mathrm{s}$. Denote by $\omega$ a supposed desired melt pool
\[\omega(t) = \{\bm{x} \in \Omega: u(\bm{x}, t) \geq u_\mathrm{melt}\}.\] 
\noindent It is difficult to express $\omega(t)$ explicitly. Instead, we consider the final solidified volume. With our purely thermal model, the powder-solid interface of this volume is dependent on the maximum temperature and would be easier to explicitly define than the melt pool $\omega(t)$. However, even further reductions can be made by isolating particular curves on this powder-solid interface. More precisely, we introduce a secondary path $\mathcal{C}^\mathrm{wd}$ chosen such that it lies on a desired powder-solid interface. The secondary path is related to the beam path by some function $\mathfrak{F}: \mathcal{C}^\mathrm{s} \rightarrow  \mathcal{C}^\mathrm{wd}$ and we write
\[\mathcal{C}^\mathrm{wd} = \{\mathfrak{F}(\bm{x}^\mathrm{s}(t)): t \in \mathcal{T}\}\]
\noindent and let $\bm{x}^\mathrm{wd} = \mathfrak{F}(\bm{x}^\mathrm{s})$. For example, if the beam path consists of one segment, then a simple example of a secondary path is
\[
\mathcal{C}^\mathrm{wd} = \{(x^\mathrm{s}(t) + w \sin \theta_1, \, y^\mathrm{s}(t) - w \cos \theta_1,\,  -d): t \in \mathcal{T}\}.
\]
\noindent The idea is that the secondary path relates to $\mathcal{C}^\mathrm{s}$ via a width $w$ and a depth $d$. With this, the description of the optimal melting reduces to two paths; $\mathcal{C}^\mathrm{s}$ (preset) and $\mathcal{C}^\mathrm{wd}$ (chosen with respect to $\mathcal{C}^\mathrm{s}$). This reductive approach is illustrated in Figure \ref{fig:reduction}.

As we shall see in the following section, the steps taken above allow us to formulate a simple optimization problem that is efficient in the sense that we, instead of tracking some desired transient melt pool $\omega(t)$ in a volume, only track two scalar values $u_\mathrm{melt}$ and $u_\mathrm{surf}$ on paths.

\subsection{Mathematical formulation}
\label{sec:mathematical_formulation}
\noindent Following the reduction in Section \ref{sec:reduction}, we are now interested in maximum temperatures on paths running along the beam path since these temperatures determine the size of the subsequent solidified volume. Before we formulate the problem, we need the following.
\begin{definition}\label{def:hatchline}
$\mathrm{(Hatch\,\,line).}$ Given a partition as in Definition \ref{def:pw}, two segments $k$ and $k+1$, $1\leq k < N-1$, are connected if $\bm{x}_k^\mathrm{f} = \bm{x}_{k+1}^\mathrm{i}$. A sequence $\{i\}_{i=k}^{K}$ of connected segments form a hatch line if $\theta_n = \theta_m$ $\forall n, m \in [k, K]$ and $\theta_{k-1} \neq \theta_k$ and $\theta_K \neq \theta_{K+1}$.
\end{definition}
\noindent The total number of hatch lines $M$ satisfies $1\leq M \leq N$. 

Define also the maximum temperature field
\begin{align*}
\mathcal{M}(\bm{x}; \bm{d}) = \max_{t\in \mathcal{T}} \{ u\left(\bm{x}, t; \bm{d}\right)\}.
\end{align*}
\noindent We want $\mathcal{M}(\bm{x}; \bm{d}) = u_\mathrm{melt}$ for all $\bm{x}$ on $\mathcal{C}_i^\mathrm{wd}$ in order to ensure a uniform and thorough melting. Similarly, $\mathcal{M}(\bm{x}; \bm{d})$ should not exceed some maximum allowed temperature $u_\mathrm{surf}$ on $\mathcal{C}^\mathrm{s}$.

The resulting objective vector becomes
\begin{equation}\label{eq:objectives} f(\bm{d}) = \big\{f_1(\bm{d}), f_2(\bm{d})\big\},\end{equation}
\noindent where
\begin{equation}
\begin{split}
f_1(\bm{d}) &= \int_{\mathcal{C}^\mathrm{wd}} \alpha(\bm{x}^\mathrm{wd}) \cdot\big(\mathcal{M}(\bm{x}^\mathrm{wd}; \bm{d}) - u_\mathrm{melt}\big)^2 \,\mathrm{d}x,\\
f_2(\bm{d}) &= \int_{\mathcal{C}^\mathrm{s}} \alpha(\bm{x}^\mathrm{s}) \cdot\big(\mathcal{M}(\bm{x}^\mathrm{s}; \bm{d}) - u_\mathrm{surf}\big)^2 \,\mathrm{d}x.
\end{split}
\label{eq:fs}
\end{equation}
\noindent Here
\begin{equation}
    \alpha(\bm{x}) = \begin{cases}
    0 \text{\,\,\, if }\bm{x}\text{ is near the start or end of a hatch line,}\\
    1 \text{\,\,\, otherwise}
    \end{cases}
    \label{eq:alpha}
\end{equation}
\noindent is a weight that excludes intervals from the cost functional if they are very close to the start point or end point of a hatch line. This type of weight is inserted because the total heat supplied to a region near a start point, for instance, is comparatively small since the beam only moves away from it rather than passing it. As a consequence, it can be difficult to reach the reference temperatures in these regions and  if included, they deteriorate the overall performance of the optimizer. Therefore, it is better to ignore these intervals in the goal functional and instead let them be covered by the contouring stage, in which the beam scans along the boundary of the shape being melted. The contouring stage also improves the surface finish of the part \cite{smith}. Effectively, this choice of $\alpha$ simply means that the domains of integration in \eqref{eq:fs} become slightly smaller.

The functionals \eqref{eq:fs} are of tracking type where we use the $L_2$ norm to minimize the distances between the actual maximum temperature and the desired maximum temperatures. The reason for tracking the surface temperature rather than penalizing only too high temperatures is that it helps restricting the shape of the melt pool. Without this restriction, one could potentially end up with an extremely wide and shallow type of melting, for instance. 

The analytic solution presented in Section \ref{sec:thermal_model} allows for pointwise computations of temperatures, and that is why we can easily compute temperatures on lines, which saves a large amount of computation time compared to doing so on surfaces or in volumes. This remark highlights a big motivation behind the reduction carried out in Section \ref{sec:reduction}.

The objective vector \eqref{eq:objectives} needs to be translated into a scalar valued cost functional in order to use standard nonlinear programming solvers. We use a scalarization known as the weighting method. In this method the weighted sum of the objectives is minimized. We introduce weights $W_i \geq 0, \,\, i = 1, 2$. The scalarized optimization problem becomes
\begin{equation}
\begin{split}
\begin{aligned} &\mathrm{minimize} \,\,\,\, J(\bm{d}) = W_1 f_1(\bm{d}) + W_2 f_2(\bm{d}) \\
&\mathrm{subject \,\, to} 
\end{aligned}
\\
\begin{aligned}
&\quad\mathrm{state\,\,eq.}\,\, \eqref{eq:pde}, &\mathrm{(PDE \, constraint)}\\
&\quad\bm{d}_\mathrm{min} \leq \bm{d} \leq \bm{d}_\mathrm{max}. &\mathrm{(Parameter \, constraint)}
\end{aligned}
\end{split}
\label{eq:wp}
\end{equation}
\noindent The choice of weights should represent the relative importance of the objectives; important objectives are weighted more heavily. 


\section{A first greedy algorithm for solving the scalarized optimization problem}
\label{sec:greedy_algorithm}
\noindent In order to speed up the optimization, we propose a method that makes use of the fact that the melt pool, and hence maximum temperatures, are localized to the beam. The proceeding involves a division of $\mathcal{T}$ into subintervals on the form $\bigcup_{i=p}^q (t_i^\mathrm{i}, t_i^\mathrm{f}]$. Local optimization problems are solved on these subintervals and optimal parameter pairs $(\sigma_p, v_p)$ are frozen sequentially. When given parameter pair(s) has been frozen, the local problem is translated in time (and space) and the initial condition is updated. As such, this greedy type of algorithm divides the optimization problem \eqref{eq:wp} into several smaller optimization problems that are faster to solve. 

The goal functional in \eqref{eq:wp} involves maximum temperatures over time near the beam path, which is a property that is local to the beam itself. Given a point $P$ on, say, $\mathcal{C}^\mathrm{s}$, it is known that $P$ will obtain its largest temperature during a time window when the beam, and the melt pool it generates, passes $P$. Therefore, it is the values of the beam parameters during this particular time window that has the highest influence on maximum temperature at $P$. The reasoning is similar for a point on $\mathcal{C}^\mathrm{wd}$, the only difference being that it takes slightly longer to reach the maximum temperature on $\mathcal{C}^\mathrm{wd}$ since heat diffusion is not an instantaneous process. Therefore, we decide to split the optimization problem \eqref{eq:wp} into multiple subproblems that are solved sequentially in time while parameter pairs are frozen as we go along.

In order to formalize the method, we make the following definition.
\begin{definition}\label{def:window}
$\mathrm{(Window).}$ Given a partition as in Definition \ref{def:pw}, a time window $\mathcal{T}_{p, q}=\cup_{i = p}^q (t_i^\mathrm{i}, t_i^\mathrm{f}]$ is defined as a set of adjacent segments in $\mathcal{T}$. The size of $\mathcal{T}_{p, q}$ is the number of segments that constitutes it. The local beam path and local secondary path corresponding to $\mathcal{T}_{p, q}$ are given by 
\begin{align*}
\mathcal{C}_{p, q}^\mathrm{s} &= \{\bm{x}^\mathrm{s}(t): t \in \mathcal{T}_{p, q}\},\\
\mathcal{C}_{p, q}^\mathrm{wd} &= \{\bm{x}^\mathrm{wd}(t): t \in \mathcal{T}_{p, q}\}.
\end{align*}
\end{definition}
\noindent Hence the beam traverses the path $\mathcal{C}_{p, q}^\mathrm{s}$ during time $\mathcal{T}_{p, q}$.

\begin{figure}[h!]
\captionsetup{width=0.9\linewidth}
\includegraphics[width=0.9\linewidth]{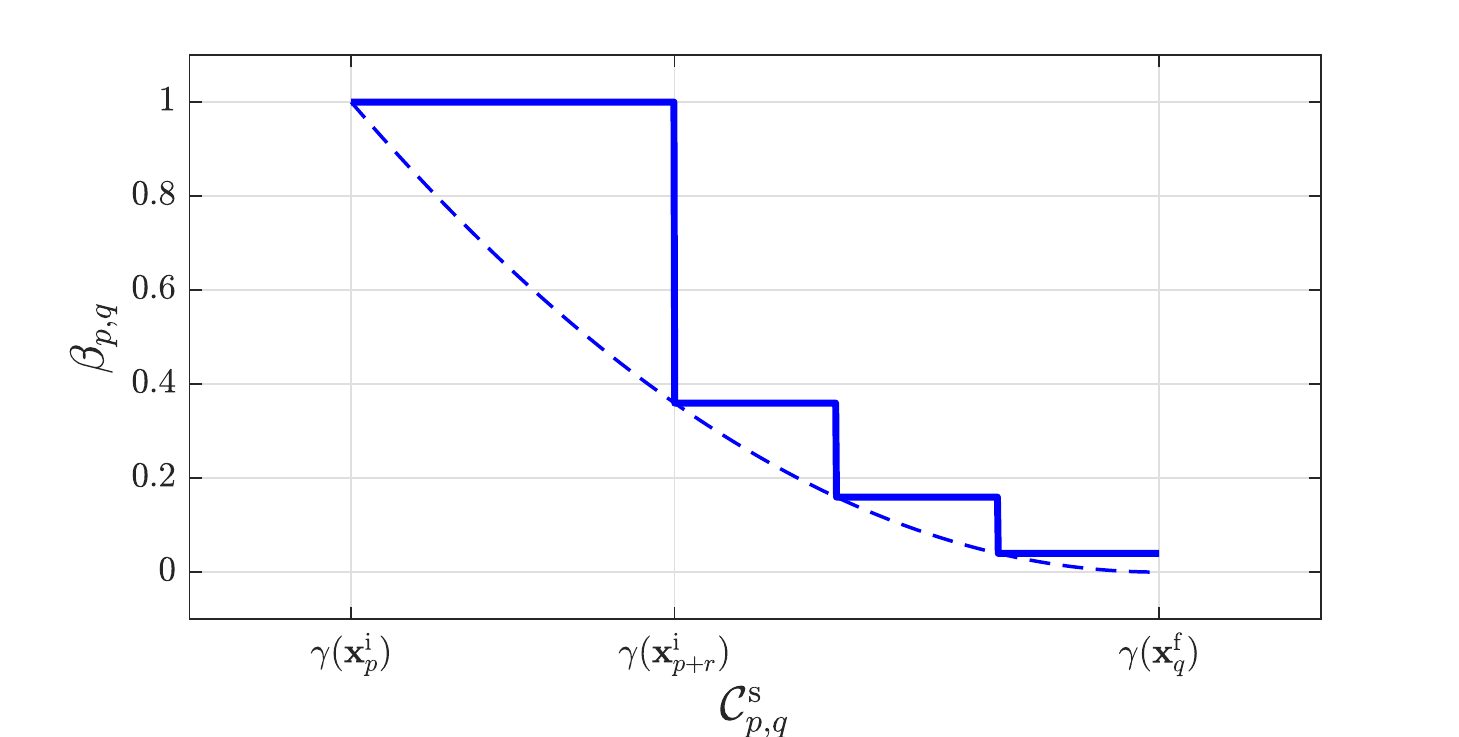}
\caption{The weight $\beta_{p, q}$ is based on a quadratic function that decreases along the local beam path $\mathcal{C}_{p,q}^\mathrm{s}$. It is piecewise constant, largest on the first $r$ segments in the window since they are about to be frozen (see Algorithm~\ref{alg:greedy_alg}), and takes different values over the remaining segments. The weight is identical for the local secondary path $\mathcal{C}_{p,q}^\mathrm{wd}$.}
\label{fig:weight}
\end{figure}

Define a local decision vector $\bm{d}_{p, q} = \{(\sigma_k, v_k)\}_{k=p}^q$ and a local maximum temperature field 
\begin{align*}
\mathcal{M}_{p,q}(\bm{x};\bm{d}_{p, q}) = \max_{t\in \mathcal{T}_{p, q}} \{ u\left(\bm{x}, t; \bm{d}_{p, q}\right)\}.
\end{align*}
\noindent Similarly, we define the local objectives
\begin{align*}
g_1(\bm{d}_{p, q}) &= \int_{\mathcal{C}_{p, q}^\mathrm{wd}} \beta_{p, q}(\bm{x}^\mathrm{wd}) \cdot \alpha(\bm{x}^\mathrm{wd}) \cdot \left(\mathcal{M}_{p, q}(\bm{x}^\mathrm{wd}; \bm{d}_{p, q}) - u_\mathrm{melt}\right)^2 \,\mathrm{d}s,\\
g_2(\bm{d}_{p, q}) &= \int_{\mathcal{C}_{p, q}^\mathrm{s}} \beta_{p, q}(\bm{x}^\mathrm{s}) \cdot \alpha(\bm{x}^\mathrm{s}) \cdot \left(\mathcal{M}_{p, q}(\bm{x}^\mathrm{s}; \bm{d}_{p, q}) - u_\mathrm{surf}\right)^2 \,\mathrm{d}s.\\
\end{align*}
\noindent Here $\beta_{p, q}$ is used to prioritize minimization of the errors over the earlier segments in the window. This weight accentuates the error on the segment(s) that is about to become frozen and it plays a crucial role; since the greedy algorithm never returns to a segment once it has been frozen it is important that the solver prioritizes this segment. Here we let $\beta_{p, q}$ be piecewise constant over the segments and determined by a function that decreases quadratically along $\mathcal{C}_{p, q}^\mathrm{s}$. See Figure \ref{fig:weight}. Formally, we have
\[\beta_{p, q}(\bm{x}^\mathrm{s}) = \chi_{(\gamma(\bm{x}_p^\mathrm{i}), \gamma(\bm{x}_{p+r}^\mathrm{i})]} + \sum_{k=p+r}^{q-1}\left(\frac{\gamma(\bm{x}^\mathrm{f}_q) - \gamma(\bm{x}_k^\mathrm{i})}{\gamma(\bm{x}^\mathrm{f}_q)  -\gamma(\bm{x}^\mathrm{k}_p)}\right)^2 \chi_{(\gamma(\bm{x}_k^\mathrm{i}), \gamma(\bm{x}_k^\mathrm{f})]}\]
\noindent on $\mathcal{C}_{p,q}^\mathrm{s}$. We define $\beta_{p, q}(\bm{x}^\mathrm{wd})$ in a similar fashion. The choice of a quadratic underlying function is based on tests that investigate how the weight affects the optimization results. 

\begin{figure}[h!]
\captionsetup{width=0.9\linewidth}
\includegraphics[width=0.9\linewidth]{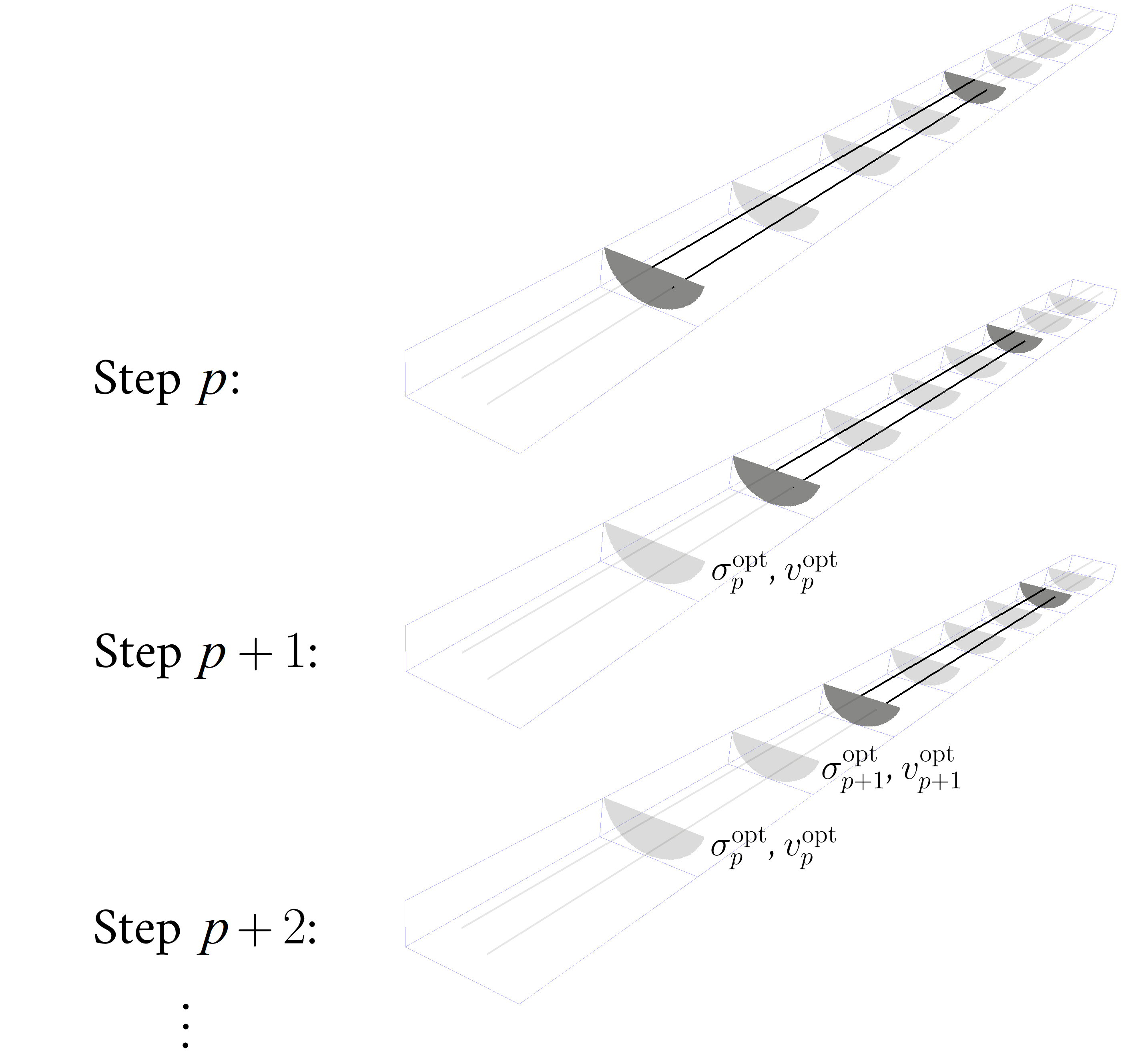}
\caption{The greedy algorithm divides the optimization problem \eqref{eq:wp} into several smaller optimization problems that are easier to solve. The procedure involves a division of the beam path into subintervals. Localized optimization problems are solved on these intervals and optimal parameters are frozen in steps. After each such step, a new sub-problem is created by translating in time (ans space) and updating the initial condition. Here $q=p+3$ and $r=1$ (see Algorithm \ref{alg:greedy_alg}).}
\label{fig:descr_greedyv2_all_num}
\end{figure}

Now, by employing the same scalarization as for the global problem \eqref{eq:wp}, the resulting scalarized subproblem becomes
\begin{equation}
\begin{split}
\begin{aligned} &\mathrm{minimize} \,\,\,\, J_{p, q}(\bm{d}_{p, q}) = W_1 g_1(\bm{d}_{p, q}) + W_2 g_2(\bm{d}_{p, q})\\
&\mathrm{subject \,\, to}
\end{aligned}
\\
\begin{aligned}
&\quad\mathrm{state\,\,eq.}\,\, \eqref{eq:pde}, &\mathrm{(PDE \, constraint)}\\
&\quad\bm{d}_{{p, q}_\mathrm{min}} \leq \bm{d}_{p, q} \leq \bm{d}_{{p, q}_\mathrm{max}}. &\mathrm{(Parameter \, constraint)}
\end{aligned}
\end{split}
\label{eq:wsp}
\end{equation}
\noindent We can now formulate the greedy algorithm. This is done in Algorithm \ref{alg:greedy_alg}, and some complementary comments are given below. An illustration of the main idea is given in Figure \ref{fig:descr_greedyv2_all_num}. 
\begin{tcolorbox}[blanker,float=btp, grow to left by=1cm,grow to right by=1cm]
\begin{algorithm}[H]
 \KwIn{A partition as in Definition \ref{def:pw}.
\newline Remaining optimization problem data ($\mathcal{C}^\mathrm{wd}, u_\mathrm{melt}, u_\mathrm{surf}$, bounds, weights).
}
\BlankLine
\Parameter{$(\bm{\sigma}, \bm{v}) = \{(\sigma_k, v_k)\}_{k=1}^N$ \Comment{Parameter pairs/decision array.}
\newline $p = 1$ \Comment{Index of first segment in current window.}
\newline $q \in \{1, \hdots, N\}$ \Comment{Index of last segment in current window.}
\newline $r \in \{1, \hdots, q-p+1\}$ \Comment{Number of segments  to freeze next.}
}
\BlankLine
\KwOut{$(\bm{\sigma}^\mathrm{opt}, \bm{v}^\mathrm{opt}), \,\, J^\mathrm{opt}$ \Comment{Optimal\footnote{With respect to the subproblems \eqref{eq:wsp}. We can not expect to find the optimal decision vector for the global problem \eqref{eq:wp}, only an approximation of it.} decision vector and objective.}
}
\BlankLine
\Begin{
\While{$p \leq N$}{
\vspace{0.1cm}
Solve subproblem \eqref{eq:wsp} for window $\mathcal{T}_{p, q}$ with initial guess 
\newline $\bm{d}_{p,q} =\{(\sigma_k, v_k)\}_{k=p}^q$ to find candidate decision variables $\{(\tilde{\sigma}_k, \tilde{v}_k)\}_{k=p}^q$.\\
\vspace{0.2cm}
 $  (\sigma_k, v_k) \leftarrow
\begin{cases}
(\tilde{\sigma}_k, \tilde{v}_k), \,\, k = p, \,\hdots, q   \\
(\tilde{\sigma}_q, \tilde{v}_q), \,\, k =q+1, \,\hdots, N \\
\end{cases}
$ \Comment{Update parameter pairs.}\\
\vspace{0.2cm}
$(\sigma_k^\mathrm{opt}, v_k^\mathrm{opt}) =  (\sigma_k, v_k), \,\, k = p, \hdots, p+r-1$\Comment{Freeze $r$ parameter pairs.}\\
\vspace{0.4cm}
Change window location:\\
$p \leftarrow p + r$\Comment{Update first index.}\\
\vspace{0.1cm}
$q \leftarrow \min\{q(p) + r, N\}$\Comment{Update last index.}\\
\vspace{0.1cm}
$r \leftarrow \min\{r(p), q-p+1\}$\Comment{Update number of segments to freeze next.}\\
}
\vspace{0.1cm}
$J^\mathrm{opt} = J((\bm{\sigma}^\mathrm{opt}, \bm{v}^\mathrm{opt}))$\Comment{Compute optimal objective in \eqref{eq:wp}.}
}
\caption{Greedy algorithm for finding an approximate solution of the scalarized problem \eqref{eq:wp}. Note that $q-p+1$ equals the size of the current window.}\label{alg:greedy_alg}
\end{algorithm}
\end{tcolorbox}
\begin{enumerate}[labelindent=0pt,itemindent=1em]
\item[$\mathscr{L}4$:] Future parameter pairs are updated as well, because the values found in the current window are likely a better guess than the initial one.
\item[$\mathscr{L}5$:] The segments corresponding to the frozen parameter pairs are removed from the window. In the implementation, certain checks can be made to determine whether freezing should take place or not (as in, the amount of pairs to freeze). We leave out the details.
\item[$\mathscr{L}7$-$9$:] The window size is updated in preparation for the next iteration. The min operator is used to handle the ending when $q=N$. Note also that the parameters $q$ and $r$ may depend on $p$ (i.e., on the location of the window). For instance, in a region where the beam path is complex or where the lengths of the segments are small, we might require a large window size $q-p+1$, and hence a large $q$. Furthermore, the size of the melt pool should be taken into account when choosing the value of $q$.
\end{enumerate}

\noindent The presented greedy algorithm can be utilized as a standalone tool for process optimization. If the beam path consists of $N$ segments, the total number of parameters to optimize becomes $2N$. Now, depending on the design of the layer being melted, the value of $N$ may be very high. In the next section we present a second version of the greedy algorithm that has the ability to significantly lower the size of the decision vector.


\section{A second greedy algorithm based on fitting beam parameter functions}
\label{sec:alt_greedy_algorithm}

The greedy algorithm of this section expands on an example in \cite{forslund}. As we shall see, it is similar to Algorithm \ref{alg:greedy_alg} in most regards, but it is based on educated guesses of how the beam parameters should behave.

Recall from \eqref{eq:parameter_functions} the expression for piecewise constant beam parameters. The subsequent optimization makes no assumption on the behavior of the beam parameters along the beam path. An alternative approach is to do curve fitting of prespecified beam parameter functions. To this end, we write
\begin{equation}
\begin{split}
\sigma(\bm{x}^\mathrm{s}) = \sum_{k=1}^{N} F^\sigma(\bm{x}_k^\mathrm{i}; \bm{\Lambda}^\sigma) \chi_{(\gamma(\bm{x}_k^\mathrm{i}), \gamma(\bm{x}_k^\mathrm{f})]}, \\
v(\bm{x}^\mathrm{s}) = \sum_{k=1}^{N} F^v(\bm{x}_k^\mathrm{i}; \bm{\Lambda}^v) \chi_{(\gamma(\bm{x}_k^\mathrm{i}), \gamma(\bm{x}_k^\mathrm{f})]},
\end{split}
\label{eq:parameter_functions_alt}
\end{equation}
\noindent where $\bm{\Lambda} = (\bm{\Lambda}^\sigma, \bm{\Lambda}^v)$, the coefficients in our beam parameter functions, become our new decision vector that we want to optimize. Given a decision vector, the beam parameters are then evaluated as (see \eqref{eq:parameter_functions})
\begin{align*}
\sigma_k &= F^\sigma(\bm{x}_k^\mathrm{i}; \bm{\Lambda}^\sigma),\\
v_k &= F^v(\bm{x}_k^\mathrm{i}; \bm{\Lambda}^v).
\end{align*}
\noindent These parameter functions are defined with respect to the hatch lines as
\[
\begin{rcases*}
F^\sigma(\bm{x}^\mathrm{s}; \bm{\Lambda}^\sigma) = F_l^\sigma(\bm{x}^\mathrm{s}; \bm{\Lambda}_l^\sigma)\\
F^v(\bm{x}^\mathrm{s}; \bm{\Lambda}^v) = F_l^v(\bm{x}^\mathrm{s}; \bm{\Lambda}_l^v)\\
\end{rcases*}
\,\,\mathrm{if}\,\,\bm{x}^\mathrm{s}\,\,\mathrm{on\,\,hatch\,\,line}\,\,l
\]
\noindent The reason for splitting $F^\sigma$ and $F^v$ between hatch lines is that a new hatch line often requires a rapid jump in beam parameter values. 

It follows from Definition \ref{def:hatchline} that hatch lines can simply be seen as an intermediate level between the segments and the beam path. Let $\mathpzc{S}:\{1, \hdots, M\} \rightarrow \{1, \hdots, N\}$ be an injective function that, given a hatch line $m$, returns the first segment in $m$. Let $\mathpzc{L}: \{1, \hdots, N\} \rightarrow \{1, \hdots, M\}$ be a surjective function that, given segment $n$, returns the hatch line that contains it. With these two functions it is possible to seamlessly work with both hatch lines and segments. For instance, the local beam path that consists of hatch lines $1$ to $3$ is $\mathcal{C}_{\mathpzc{S}(1),\,\mathpzc{S}(4)-1}^\mathrm{s}$.

\begin{figure}[h!]
\captionsetup{width=0.9\linewidth}
\includegraphics[width=0.9\linewidth, trim={0 0 6.5cm 0}]{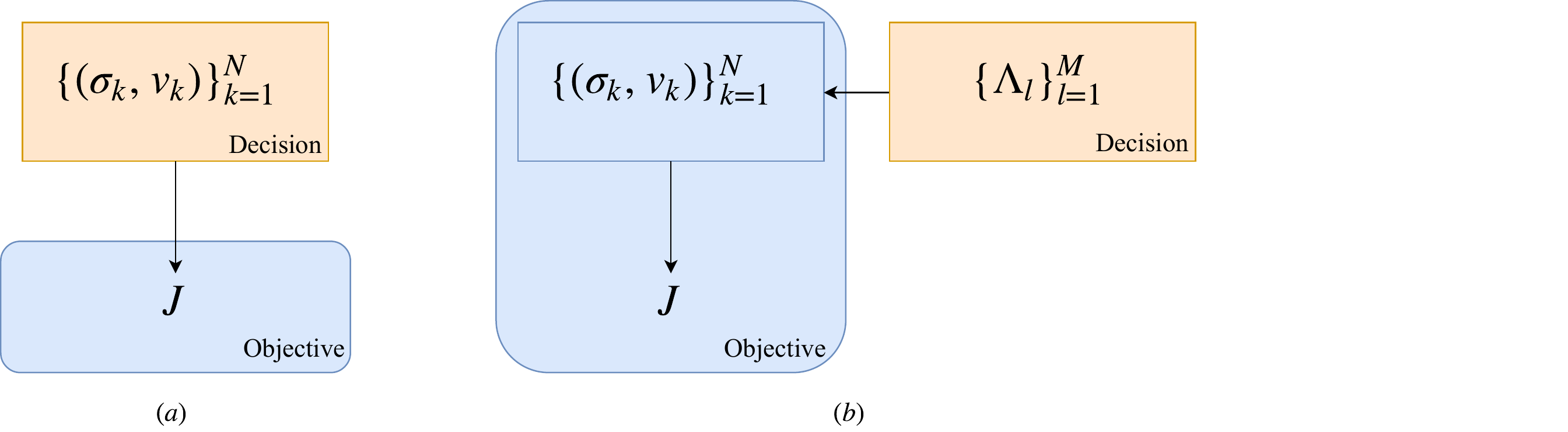}
\caption{Difference between the two greedy algorithms. In the second version (b), the additional step significantly reduces the dimension of the decision vector if the number of lines is much smaller than the number of segments, i.e., if $M \ll N$.}
\label{fig:algorithm_comparison}
\end{figure}

In terms of implementation, the second greedy algorithm is in many ways similar to the first greedy algorithm from the previous Section  \ref{sec:greedy_algorithm}. They both rely on Definition \ref{def:pw} and piecewise constant beam parameters and they both solve subproblems of the form \eqref{eq:wsp}. What separates them is the content of the decision vector $\bm{d}$ as illustrated in Figure \ref{fig:algorithm_comparison}.  In essence, the first algorithm optimizes the beam parameters \emph{segment wise} while the second algorithm optimizes the beam parameters \emph{line wise}. The second greedy algorithm is detailed in Algorithm \ref{alg:greedy_alg_alt}.
\begin{tcolorbox}[blanker,float=btp, grow to left by=1cm,grow to right by=1cm]
\begin{algorithm}[H]
 \KwIn{A partition as in Definition \ref{def:pw}; $M$ hatch lines.
\newline Remaining optimization problem data ($\mathcal{C}^\mathrm{wd}, u_\mathrm{melt}, u_\mathrm{surf}$, bounds, weights).
}
\BlankLine
\Parameter{$\bm{\Lambda} = \{\bm{\Lambda}_l\}_{l=1}^M$ \Comment{Coefficient sets/decision array.}
\newline $p = 1$ \Comment{Index of first hatch line in current window.}
\newline $q \in \{1, \hdots, M\}$ \Comment{Index of last hatch line in current window.}
\newline $r \in \{1, \hdots, q-p+1\}$ \Comment{Number of hatch lines to freeze next.}

}
\BlankLine
\KwOut{$\bm{\Lambda}^\mathrm{opt}, \,\, J^\mathrm{opt}$ \Comment{Optimal\footnote{With respect to the subproblems \eqref{eq:wsp}. We can not expect to find the optimal decision vector for the global problem \eqref{eq:wp}, only an approximation of it.} decision vector and objective.}
}
\BlankLine

\Begin{
\While{$p \leq N$}{
\vspace{0.1cm}
Solve subproblem \eqref{eq:wsp} for window $\mathcal{T}_{\mathcal{S}(p),\,\mathcal{S}(q+1)-1}$ with initial guess 
\newline $\bm{d}_{\mathcal{S}(p),\,\mathcal{S}(q+1)-1} =\{\bm{\Lambda}_l\}_{l=p}^q$ to find candidate decision variables $\{\tilde{\bm{\Lambda}}_l\}_{l=p}^q$.\\
\vspace{0.2cm}
 $ \bm{\Lambda}_l \leftarrow
\begin{cases}
\tilde{\bm{\Lambda}}_l, \,\, l = p, \,\hdots, q   \\
\tilde{\bm{\Lambda}}_q, \,\, l =q+1, \,\hdots, M \\
\end{cases}
$ \Comment{Update parameter pairs.}\\
\vspace{0.2cm}
$\bm{\Lambda}_l^\mathrm{opt} = \bm{\Lambda}_l, \,\, l = p, \hdots, p+r-1$\Comment{Freeze $r$ parameter pairs.}\\
\vspace{0.4cm}
Change window location:\\
$p \leftarrow p + r$\Comment{Update first index.}\\
\vspace{0.1cm}
$q \leftarrow \min\{q(p) + r, N\}$\Comment{Update last index.}\\
\vspace{0.1cm}
$r \leftarrow \min\{r(p), q-p+1\}$\Comment{Update number of segments to freeze next.}\\
}
\vspace{0.1cm}
$J^\mathrm{opt} = J(\bm{\Lambda}^\mathrm{opt})$\Comment{Compute optimal objective in \eqref{eq:wp}.}
}
\caption{Second greedy algorithm for finding an approximate solution of the scalarized problem \eqref{eq:wp}. It is in many ways similar to the first Algorithm \ref{alg:greedy_alg}, but uses a different decision vector $\bm{d}$. Note that indices $p$, $q$ and $r$ now count over hatch lines instead of over segments.}
\label{alg:greedy_alg_alt}
\end{algorithm}
\end{tcolorbox}


\section{Numerical examples and discussion}
\label{sec:examples}
\noindent We apply the greedy algorithm on a couple of single layer problems. The scalarized subproblems \eqref{eq:wsp} are solved with the L-BFGS-B optimization algorithm \cite{Zhu:1997:ALF:279232.279236, doi:10.1137/0916069} provided by the optimization package of SciPy \cite{scipy}. Given some iterate $\hat{\bm{d}}_{p, q}$, \texttt{scipy.optimize} approximates the gradient of $J_{p,q}(\hat{\bm{d}}_{p, 1})$ using a 2-point finite difference estimation. Then L-BFGS-B, which is a quasi-Newton method, approximates the Hessian that enters in the local quadratic approximation of $J_{p,q}(\hat{\bm{d}}_{p, q})$. The solver options are selected to fit the scale of the problems considered here. 

It is worth noting that the objective $J_{p, q}$ is not differentiable. At an initial stage, not only L-BFGS-B but also some gradient free methods were tested, and L-BFGS-B performed the best out of all solvers in those trials.  It is not remarkable that L-BFGS-B performs well on nonsmooth problems as well (although we can not expect the same convergence as for a smooth optimization problem) \cite{Lewis2013}. On a related note, in the implementation we relax the scalarized subproblem \eqref{eq:wsp} somewhat by replacing the (local) maximum temperature field $\mathcal{M}_{p,q}(\bm{x};\bm{d}_{p, q})$ with an approximation,
\[
\mathcal{M}_{p,q}(\bm{x};\bm{d}_{p, q}) \approx \frac{1}{K} \log \left( \int_{\mathcal{T}_{p, q}} \exp(K \cdot u(\bm{x}, s; \bm{d}_{p, q})) \,\mathrm{d}s \right),
\]
\noindent for an appropriate scalar $K$, which improves performance slightly. Finally, while the choice of starting point/initial decision vector can have a large impact on the performance of the optimizer, efforts related to this choice are not the main focus here and so disregarded.

It is important to emphasize that for practical use of the optimization scheme, the material data that enter the thermal model need to be fit with respect to experiments or a more detailed model. The determination of effective parameters is crucial since the model is simple and based on several assumptions. In the following examples we use material parameters that represent Ti-$6$Al-$4$V.

In the current implementation of the beam scanning, there are no pauses between any segments during melting. For instance, a jump from one hatch line to the next is instantaneous. However, adding delay time for jumps is straightforward. 


\subsection{Example 1: segment wise optimization on snake pattern}
\label{sec:example1}
\begin{figure}[h!]
\captionsetup{width=0.9\linewidth}
\includegraphics[width=0.9\linewidth]{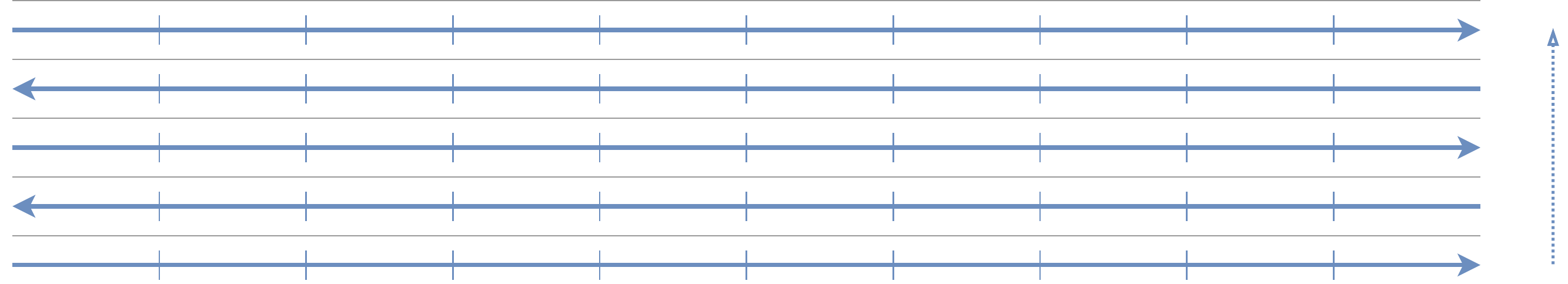}
\caption{Beam path in Example 1: 5 hatch lines with length $5$ mm. The path consists of $50$ segments that are separated by the small ticks. The line offset (i.e., the distance between two adjacent lines) is $200\, \mu\mathrm{m}$. The thin gray lines indicate the secondary path $\mathcal{C}^\mathrm{wd}$.}
\label{fig:ex1_beampath}
\end{figure}

\noindent We illustrate how the optimization scheme resolves certain heating related issues. The beam path is shown in Figure \ref{fig:ex1_beampath}. It consists of $50$ segments, each with length $0.5\,\mathrm{mm}$. The line offset during hatching is $l_\mathrm{off} = 200\, \mu\mathrm{m}$. The hatching is performed in a snake-like manner in the upward $y$-direction. The weight $\alpha(\bm{x})$ (see \eqref{eq:alpha}) is zero on the first $0.4\,\mathrm{mm}$ and last $0.4\,\mathrm{mm}$ of a hatch line. While this beam path amounts to a simple rectangular shape, it still allows for several types of investigations.

\begin{figure}[h!]
\captionsetup{width=0.9\linewidth}
\includegraphics[width=0.9\linewidth]{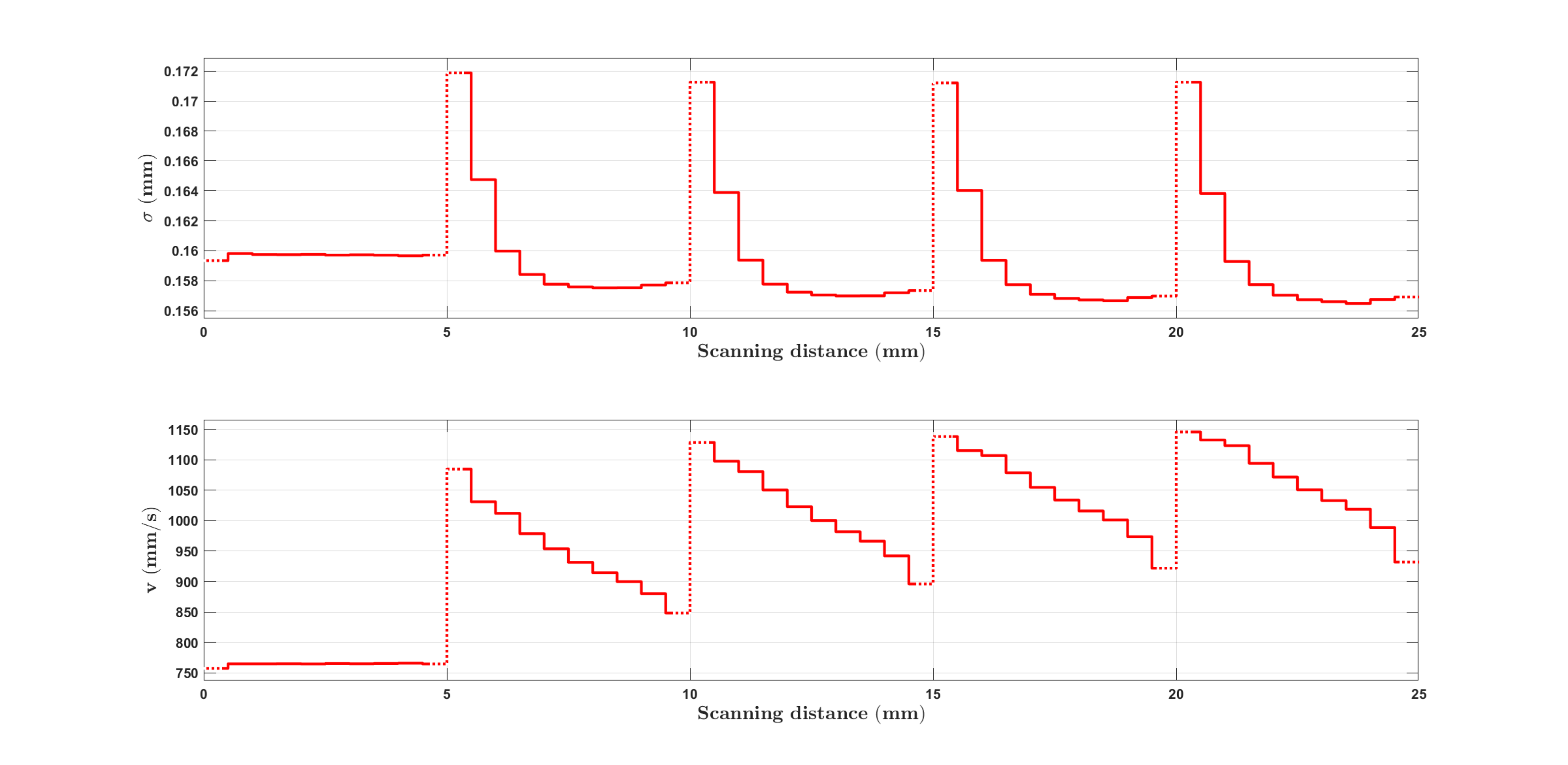}
\caption{The solution as obtained by the first greedy algorithm. The dotted parts indicate the intervals where the weight $\alpha$ is $0$ (the start and end of the hatch lines).}
\label{fig:ex1_optimal_parameters}
\end{figure}

\begin{figure}[h!]
\captionsetup{width=0.9\linewidth}
\includegraphics[width=0.9\linewidth]{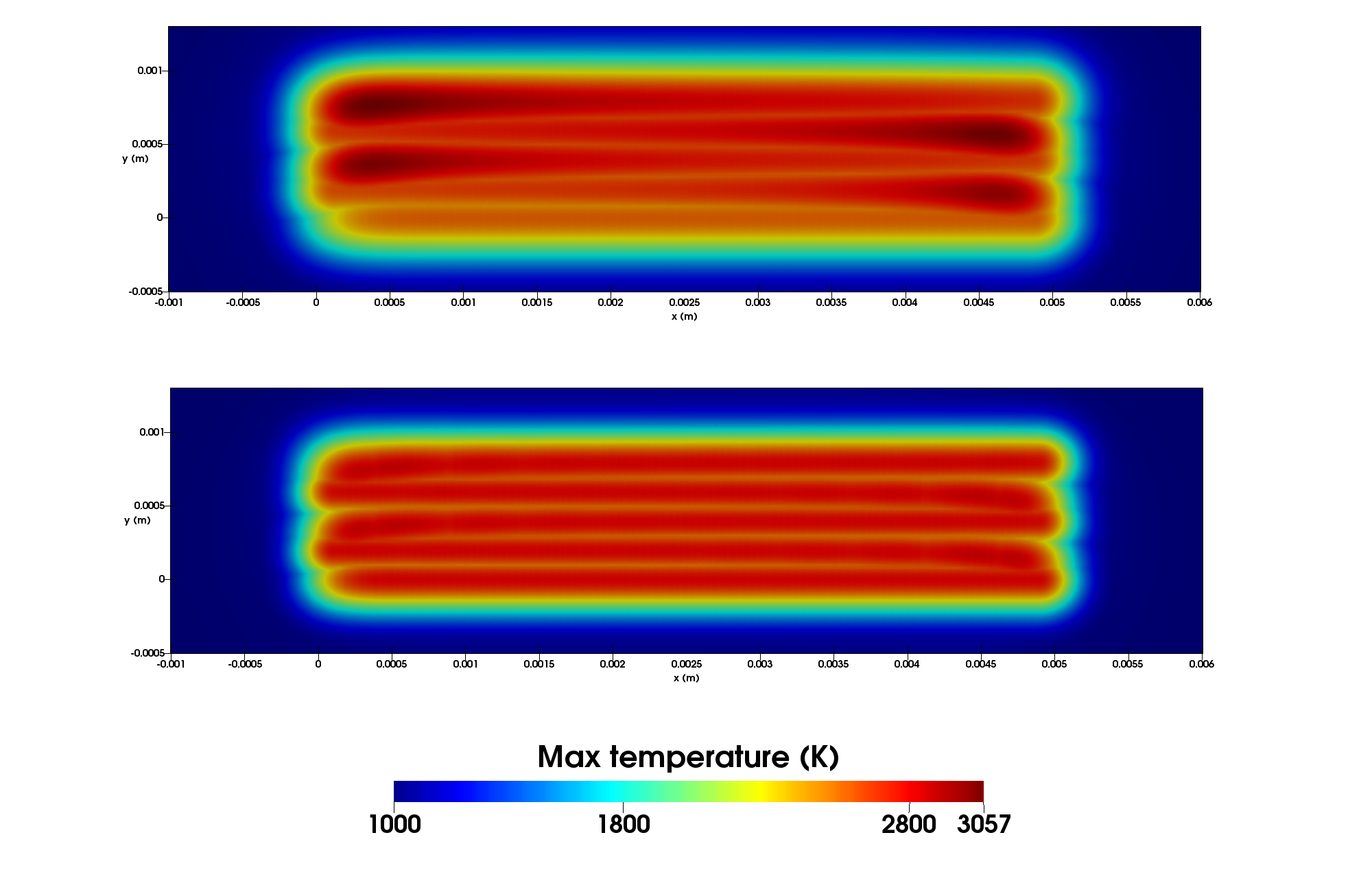}
\caption{A comparison between maximum temperatures on the surface ($z = 0$) before optimization (top) and after optimization. The optimization scheme resolves the heat accumulation at the turning points.}
\label{fig:ex1_1_all}
\end{figure}

\begin{figure}[h!]
\captionsetup{width=0.9\linewidth}
\includegraphics[width=0.9\linewidth]{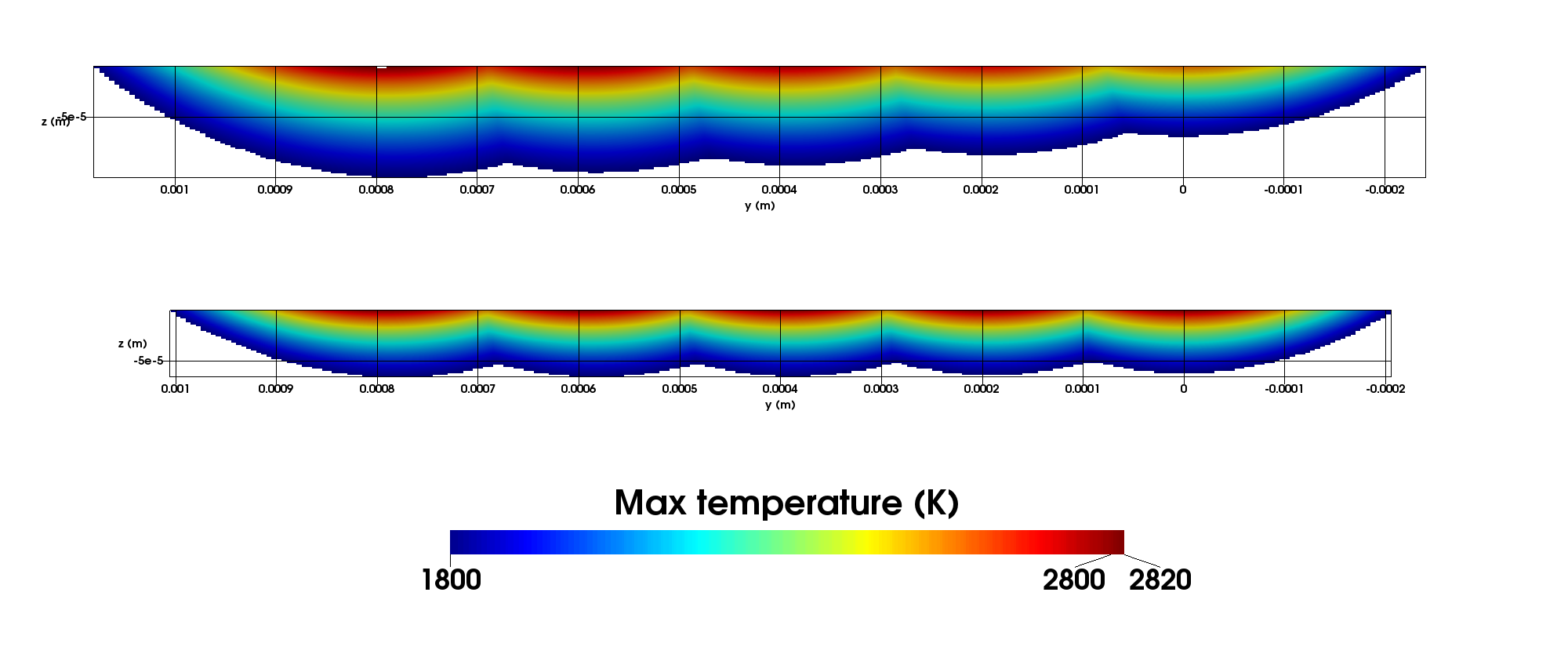}
\caption{A comparison between maximum temperatures in the cross-section $x = 2.5$ mm before optimization (top) and after optimization. Initially, the melt area increases for each hatch line since the heat influx is larger than the rate of diffusion. The optimization scheme resolves this issue and makes the area more uniform.}
\label{fig:ex1_2_all}
\end{figure}

\begin{figure}[h!]
\captionsetup{width=0.9\linewidth}
\includegraphics[width=0.9\linewidth]{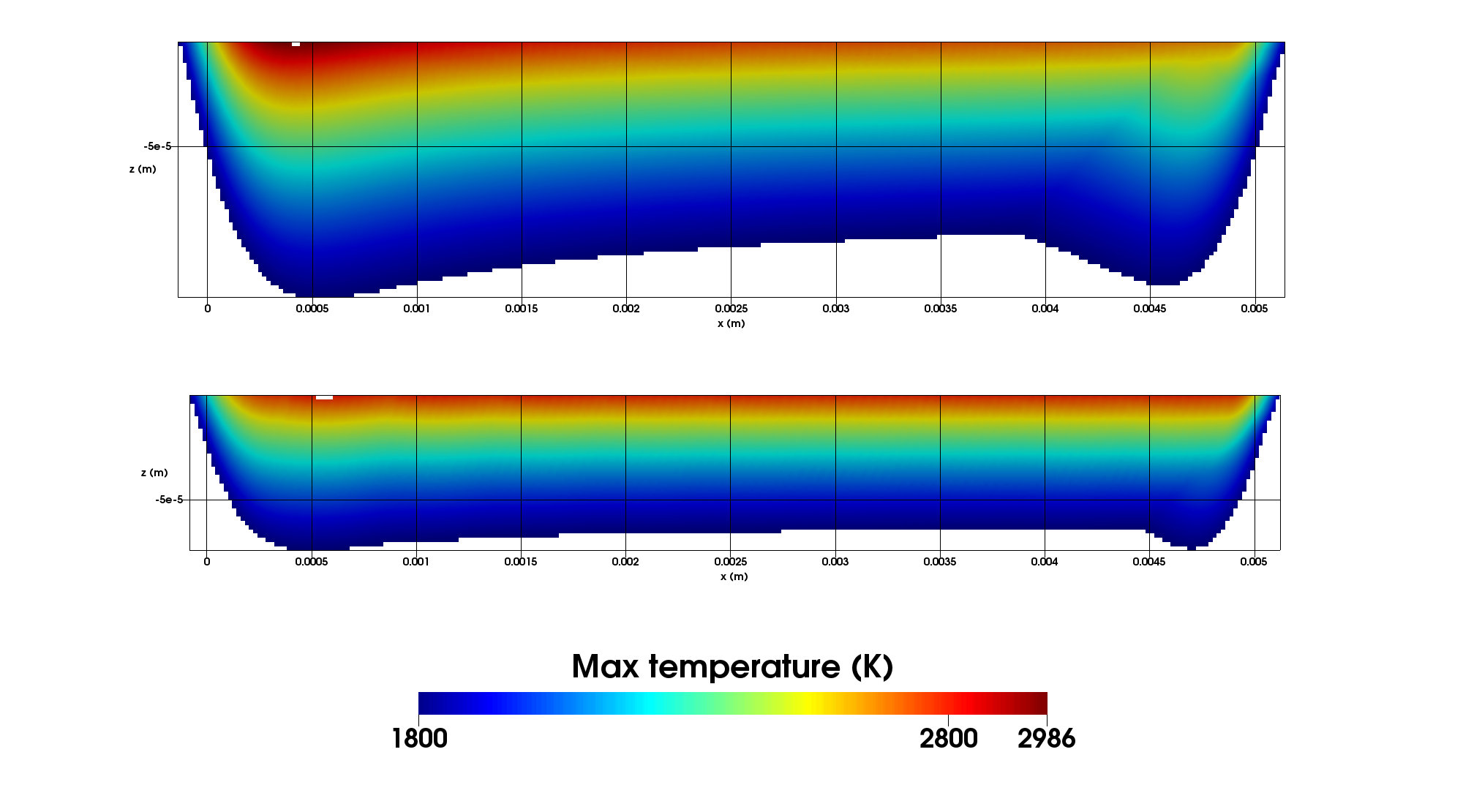}
\caption{A comparison between maximum temperatures in the cross-section $y = 2 \, l_\mathrm{off} =  0.4$ mm (i.e., along the $3^{rd}$ hatch line) before optimization (top) and after optimization. The depth of the melt area is reduced and made more uniform. The domain is scaled by a factor $10$ in the vertical direction.}
\label{fig:ex1_4_all}
\end{figure}

The secondary path is set to
\[\mathcal{C}^\mathrm{wd} = \{(x^\mathrm{s}(t), \, y^\mathrm{s}(t) + w,\,  - d): t \in \mathcal{T}\}.
\]
\noindent To avoid porosity, $w$ and $d$ need be chosen such that unmelted gaps between lines are avoided. Here we set $w = l_\mathrm{off}/2 = 100 \, \mu\mathrm{m}$ and $d = 50 \, \mu\mathrm{m}$. A complete list of parameter values is given in Table \ref{tab:ex1}.
\begin{table}[h!]
\centering
\begin{tabular}{l*{4}{l}r}
\hline
\hline
{\bf PDE specific input}   &   &   &     \\
\hline
Thermal conductivity & $\lambda$ & $20$ & $(\mathrm{W/mK})$  \\
Thermal diffusivity  & $\kappa$ & $8.45\mathrm{e}{\scalebox{1.2}[1.2]{-}}6$ & $(\mathrm{m^2/s})$  \\
Initial temperature  & $u_\mathrm{init}$ & $1000$ & $(\mathrm{K})$  \\
Absorbed beam power  & $P$ & $100$ & $(\mathrm{W})$  \\
\hline
{\bf Greedy algorithm specific input}              &  & &     \\
\hline
Reference surface temperature           & $u_\mathrm{melt}$ & $1800$ & $(\mathrm{K})$  \\
Reference melt temperature            & $u_\mathrm{surf}$ & $2800$ & $(\mathrm{K})$  \\
Secondary beam position, width           & $w$ & $100$ & $(\mu \mathrm{m})$  \\
Secondary beam position, depth     & $d$ & $50$ & $(\mu \mathrm{m})$  \\
Weight 1                         & $W_1$ & $0.7$ &   \\
Weight 2                         & $W_2$ & $0.3$ &   \\
Window size               & $q-p+1$ & $5$ &   \\
Segments frozen in each iteration               & $r$ & $1$ &   \\
Initial spot size on segment $k$           & $\sigma_k$& $0.2\,\,\forall \,k$ & $(\mathrm{mm})$ \\
Initial speed on segment $k$               & $v_k$& $0.5\,\,\forall \,k$ & $(\mathrm{m/s})$ &   \\
Bounds, spot size           & $(\sigma_\mathrm{min}, \sigma_\mathrm{max})$& $(1\mathrm{e}{\scalebox{1.2}[1.2]{-}}2, 1\mathrm{e}0)$ & $(\mathrm{mm})$ \\
 Bounds, speed               & $(v_\mathrm{min}, v_\mathrm{max})$& $(1\mathrm{e}1, 1\mathrm{e}4)$ & $(\mathrm{mm/s})$ &   \\
\hline
\hline
\end{tabular}
\caption{Parameter values in Example \ref{sec:example1}.}
\label{tab:ex1}
\end{table}
Figure \ref{fig:ex1_optimal_parameters} shows the solution of \eqref{eq:wp} as obtained by the greedy Algorithm \ref{alg:greedy_alg}. We see rapid variations in the optimized beam parameters, in particular close to the turning points where they seek to neutralize the concentrated influx of heat that occur in those regions. Problem \eqref{eq:pde} is then solved for the optimized beam parameters and the resulting maximum temperature in various slices of the domain is shown in Figures \ref{fig:ex1_1_all}, \ref{fig:ex1_2_all}, and \ref{fig:ex1_4_all}. These figures include the initial maximum temperature for comparison. 
\noindent The results indicate that despite the reductions leading up to its formulation, the greedy algorithm is able to control the heat generated during melting to a rather large degree. 

One concern with the greedy algorithm is that it carries with it several uncertainties. Many trials are required to find proper values for the parameters that make up the scheme, such as $\alpha$, the window size and segment lengths, since they depend on the beam path and thermal diffusivity. 


\subsection{Example 2: segment wise optimization on nonparallel pattern}
\label{sec:example2}
\begin{figure}[h!]
\captionsetup{width=0.9\linewidth}
\includegraphics[width=0.9\linewidth]{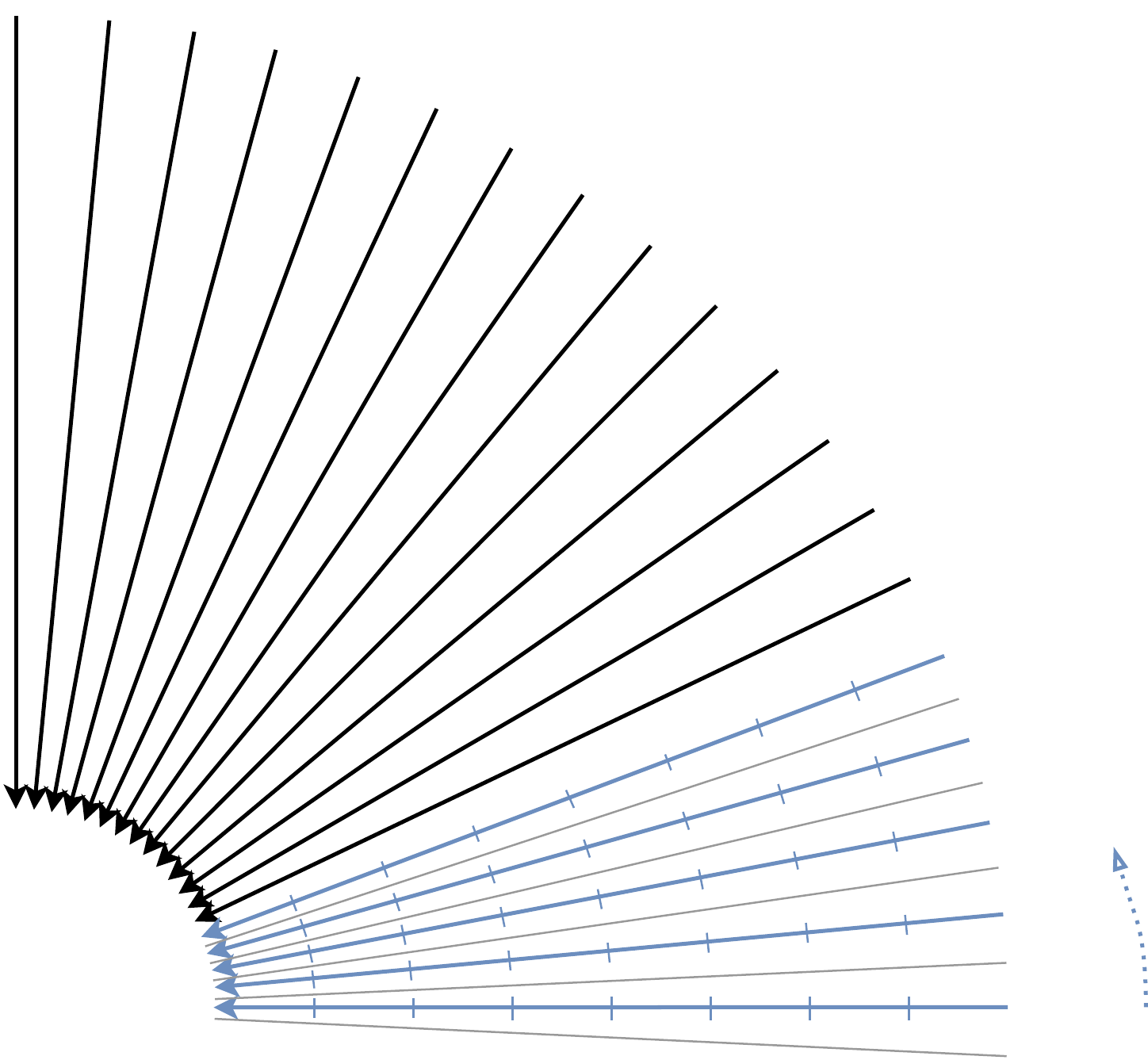}
\caption{Beam path in Example 2 for melting the first quadrant of an annulus with $r_i = 1$ mm and outer radius $r_o = 5$ mm.  The entire beam path consists of $19$ lines of length $4$ mm. Each line is divided into $8$ segments of equal length $0.5$ mm. Only the blue part, $\mathcal{C}^\mathrm{s}$, of the entire beam path is considered during optimization and the results from this optimization is then extended to the entire beam path. The path consists of $40$ segments that are separated by the small ticks. The thin gray lines indicate the secondary path $\mathcal{C}^\mathrm{wd}$.}
\label{fig:ex2_beampath}
\end{figure}

\noindent We melt the first quadrant of an annulus. The annulus has an inner radius $r_i=1\,\mathrm{mm}$ and outer radius $r_o=5\,\mathrm{mm}$. The beam path is shown in Figure \ref{fig:ex2_beampath}. 
\noindent It consists of 19 lines of length $4\,\mathrm{mm}$. Two adjacent lines differ by an angle of $5{}^\circ$. Consequently, the distance between them at $r_i$ and $r_o$ are about $87\,\mu\mathrm{m}$ and $436\,\mu\mathrm{m}$, respectively. Each line is divided into $8$ segments of equal length $0.5\,\mathrm{mm}$. The hatching is performed in the counter-clockwise direction. The thin gray lines indicate the secondary path $\mathcal{C}^\mathrm{wd}$. Once again, the weight $\alpha(\bm{x})$ is zero on the first $0.4\,\mathrm{mm}$ and last $0.4\,\mathrm{mm}$ of a hatch line. The remaining parameters are identical to the ones used in the previous example and are listed in Table \ref{tab:ex1}.

Only the blue part of the beam path is considered during optimization. More precisely, we apply the greedy algorithm on the first $5$ lines only. The results from this optimization is then extended by letting the beam parameters on lines $6$-$19$ equal the optimal beam parameters on line $5$. 

\begin{figure}[h!]
\captionsetup{width=0.9\linewidth}
\includegraphics[width=0.9\linewidth]{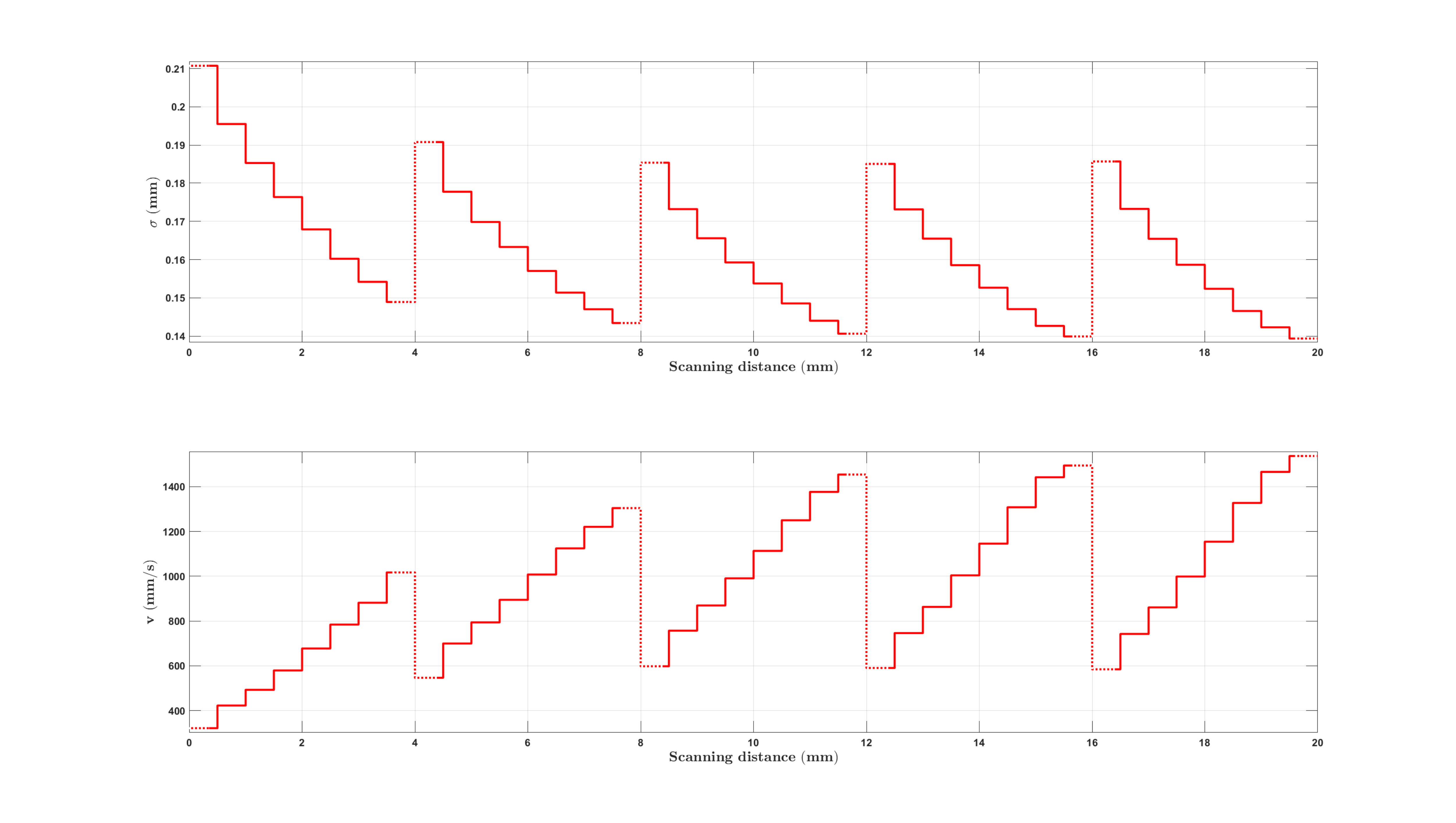}
\caption{Optimal beam parameters in the second example, as obtained by the first greedy algorithm. The dotted parts indicate the places where the weight $\alpha$ is $0$ (the start and end of the hatch lines). The parameter values start to stabilize somewhat by the fifth line. This motivates us to copy the parameter values for the fifth line to the remaining lines $6$-$19$ that were omitted in the optimization.}
\label{fig:ex2_optimal_parameters}
\end{figure}

Figure \ref{fig:ex2_optimal_parameters} shows the solution of the optimization problem. The speed increases along each hatch line since the width between a hatch line and the corresponding secondary path decreases as the hatch line approaches $r=1\,\mathrm{mm}$ (see Figure \ref{fig:ex2_beampath}).

\begin{figure}[h!]
\captionsetup{width=0.9\linewidth}
\includegraphics[width=0.9\linewidth]{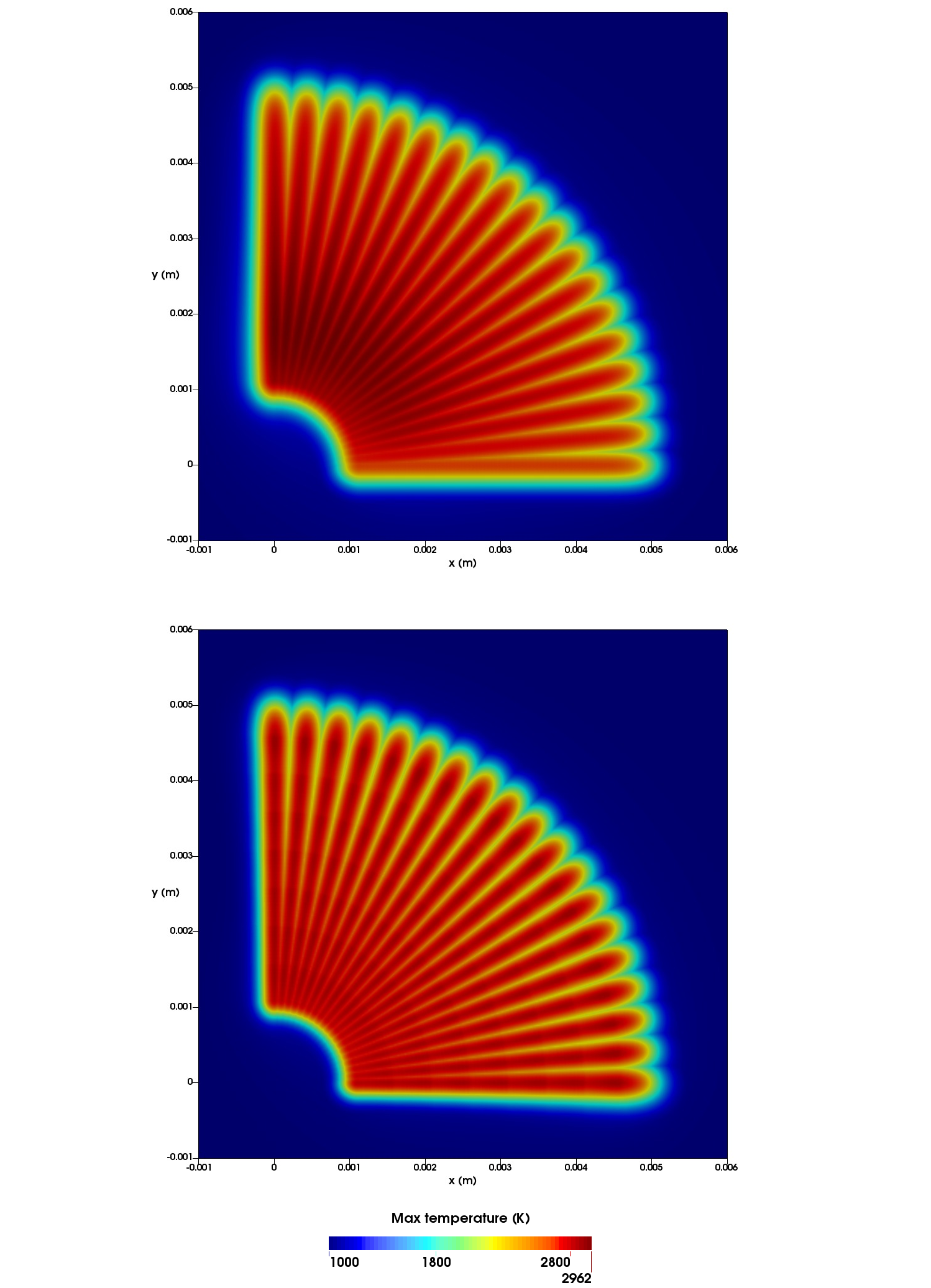}
\caption{A comparison between maximum temperatures on the surface ($z = 0$) before optimization (top) and after optimization. The optimization scheme resolves the heat accumulation near the inner radius of the annulus. The optimized maximum surface temperature appears slightly jagged along the beam path, which suggests that the segment length of $0.5$ mm is too big.}
\label{fig:ex2_all}
\end{figure}

\begin{figure}[h!]
\captionsetup{width=0.9\linewidth}
\includegraphics[width=0.9\linewidth]{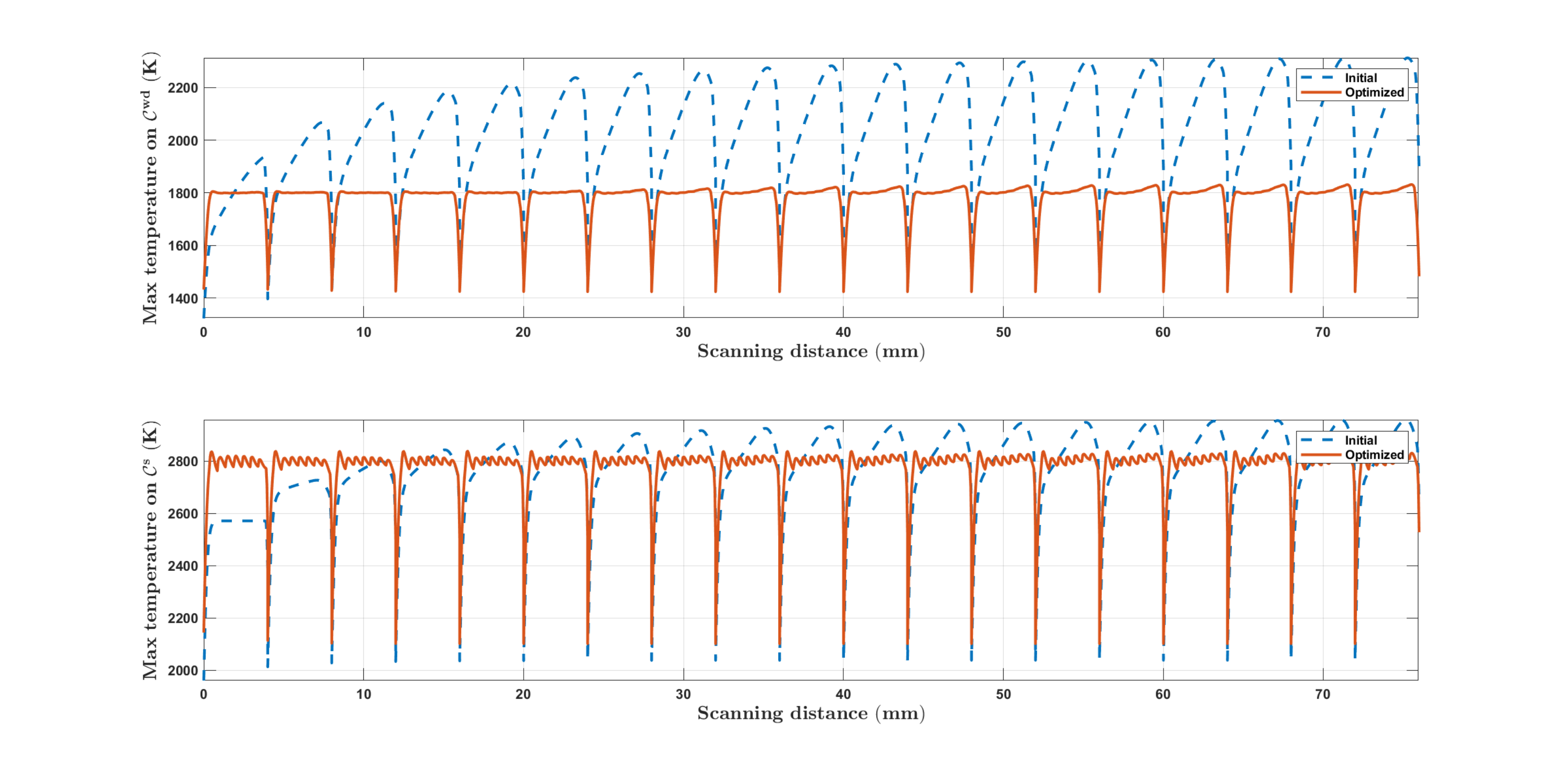}
\caption{A comparison between maximum temperatures on the (extended) secondary path and (extended) beam path. The reference temperatures $u_{\mathrm{melt}} = 1800$ K and $u_{\mathrm{surf}} = 2800$ K are successfully tracked.}
\label{fig:ex2_onpaths}
\end{figure}

After having extended this solution to the entire beam path, problem \eqref{eq:pde} is solved for the optimized beam parameters and the resulting maximum temperature is shown in Figures \ref{fig:ex2_all} and \ref{fig:ex2_onpaths}. These figures include the initial maximum temperatures for comparison. 

The plots in Figure \ref{fig:ex2_onpaths} confirm that the extension of the solution onto remaining lines $6$-$19$ works well. While the temperature on $\mathcal{C}^\mathrm{wd}$ increases near the inner radius of the annulus as the scanning progresses, this increase is small and does not justify applying the greedy algorithm on the entire beam path, $19$ lines, rather than just $5$ lines. This is just a small example of how the results from the greedy algorithm on a very small section can be utilized on larger sections of the build area. In general, this procedure offers an efficient method for improving process control: first examine and optimize typical problematic melting scenarios, then combine the results and extend them to the remainder of the layer. 


\subsection{Example 3: line wise optimization on snake pattern}
\label{sec:example3}
\begin{figure}[h!]
\captionsetup{width=0.9\linewidth}
\includegraphics[width=0.9\linewidth]{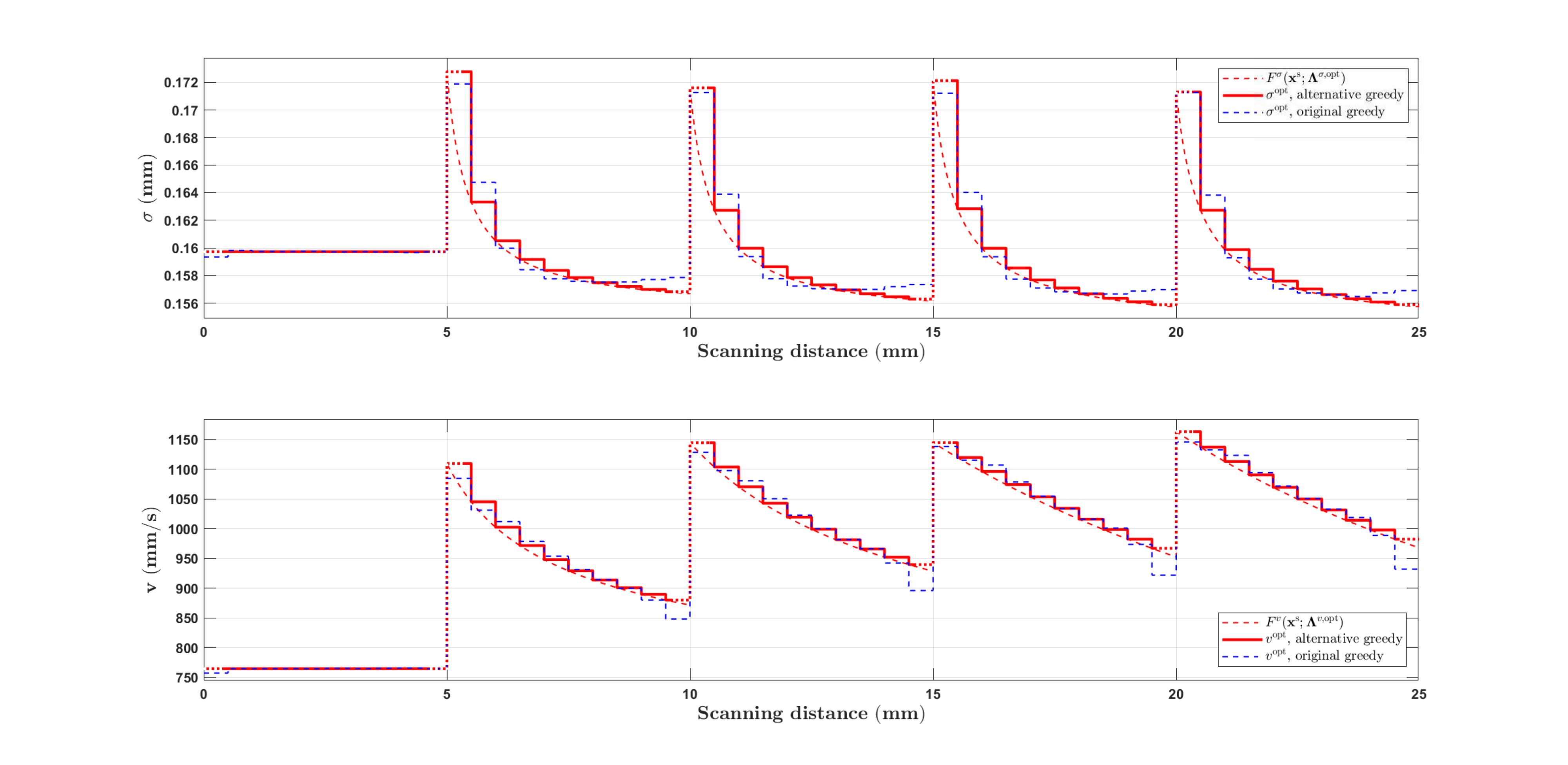}
\caption{The solution as obtained by the second greedy algorithm. The solution is a piecewise constant discretization of the beam parameter functions $F^\sigma$ and $F^v$ in \eqref{eq:preset_bpf}. The solution obtained by the first greedy algorithm is included for comparison.}
\label{fig:ex1_optimal_parameters_comparison}
\end{figure}

We now use the second greedy algorithm to solve the problem introduced in the first example, Section \ref{sec:example1}, and compare the results with the results obtained by the first greedy algorithm. Motivated by the solution obtained by the first greedy algorithm, see Figure \ref{fig:ex1_optimal_parameters}, we make the following ansatz. Let
\[\bm{\Lambda}_l = (\bm{\Lambda}_l^\sigma, \bm{\Lambda}_l^v) = (C_{1, l}^\sigma, C_{2, l}^\sigma, C_{3, l}^\sigma, C_{4, l}^\sigma, C_{1, l}^v, C_{2, l}^v, C_{3, l}^v, C_{4, l}^v)\]
\noindent and
\begin{equation}
\begin{split}
F_l^\sigma(\bm{x}^\mathrm{s}; \bm{\Lambda}_l^\sigma) &= C_{1, l}^\sigma\left(1 + \frac{C_{2, l}^\sigma}{1+C_{3, l}^\sigma(\gamma({\bm{x}^\mathrm{s}}) - \gamma(\bm{x}_{\mathcal{S}(l)}^\mathrm{i}))^{C_{4, l}^\sigma}}\right),\\
F_l^v(\bm{x}^\mathrm{s}; \bm{\Lambda}_l^v) &= C_{1, l}^v\left(1 + \frac{C_{2, l}^v}{1+C_{3, l}^v(\gamma({\bm{x}^\mathrm{s}}) - \gamma(\bm{x}_{\mathcal{S}(l)}^\mathrm{i}))^{C_{4, l}^v}}\right),
\end{split}
\label{eq:preset_bpf}
\end{equation}
\noindent for $l = 1, \hdots, M$. Hence we associate $8$ coefficients with each hatch line. From the beam path in Figure \ref{fig:ex1_beampath}, we have $M=5$. We solve optimization problem \eqref{eq:wp} according to Algorithm \ref{alg:greedy_alg_alt}. The window always consist of $1$ hatch line, i.e., $10$ segments. The results are shown in Figure \ref{fig:ex1_optimal_parameters_comparison}. The solution obtained by the first greedy algorithm is included for comparison. The results are similar. The corresponding optimal objectives $J^\mathrm{opt}$ are similar as well, as $J^\mathrm{opt}=124.18$ with the second greedy algorithm and $J^\mathrm{opt}=123.03$ with the first greedy algorithm ($J^\mathrm{init} = 1655.21$).

We end this section with a comparison between the first greedy algorithm and second greedy algorithm. The first algorithm from Section \ref{sec:greedy_algorithm} optimizes the beam parameters segment wise. It can be applied to general beam paths and there are no restrictions on the window. In particular it and can be used to look at specific problematic areas of the layer being melted in order to get an understanding of how the beam parameter functions should behave in those areas. 

The second greedy algorithm from Section \ref{sec:alt_greedy_algorithm} optimizes the beam parameters line wise. This makes the second algorithm preferable in practical problems because it significantly decreases the dimension of the decision space; the number of lines $M$ is much lower than the number of segments $N$. In the above example, the hatch lines are fairly short, but we still get $40=8M<2N=100$ when comparing dimensions of the two decision spaces. For more realistic problems the number of segments could be orders of magnitude larger, making the second greedy algorithm more attractive. 

The potential drawback of the second algorithm is that it might be difficult to make the initial ansatz for the beam parameter functions since they depend on the beam path (and secondary path). However, this is where the first algorithm can be utilized to give an initial estimate that shows the behavior of the desired beam parameters. This approach is what enabled us to choose the beam parameter functions in \eqref{eq:preset_bpf}. For the future, we imagine the development of a database of parameter functions that have been generated for different melting scenarios and that can be shared and used in other, more detailed models.


\section{Conclusions}
\label{sec:conclusions}
\noindent We have formulated an optimization scheme for controlling the heat conduction during the melting process in powder-bed-based additive manufacturing. The scheme is efficient because it exploits that the melt pool is local to the beam and shows good capabilities despite several simplifications. The current choice of objectives prioritizes speed since it only requires temperature evaluations on lines rather than in entire volumes. The scheme should be useful for studying problematic areas of the melting process where particular care needs to be put into the choice of beam parameters. 

The optimization scheme relies on a greedy algorithm. Two versions of a greedy algorithm have been presented and applied in this paper. The first one carries out optimization segment wise, which makes it applicable to many types of beam paths. The second one carries out optimization line wise, which can significantly reduce the dimension of the decision space. While the two algorithms are similar, they serve different purposes and we have detailed how they can be combined to improve process control. 

The examples considered in this paper are fairly small. For more realistic problems the number of segments could be orders of magnitude larger. By design, the greedy algorithm becomes more attractive as the total amount of segments in the beam path increases; the division of the global problem \eqref{eq:wp} into subproblems \eqref{eq:wsp} becomes more beneficial, relatively speaking, as $N$ (and $M$) increases. Furthermore, the examples are purely numerical. It is currently difficult to compare the results in Section \ref{sec:examples} with experiments because existing machines lack the functionality required to match our numerical results. Because of this, a crucial next step is to implement the necessary code in the machine such that experimental validation becomes possible. Experiments are also needed for the generation of effective parameters; the analytic model is very simple and since it contains few parameters, they need to be carefully fit against experiments.

The optimization method also needs to be complemented with different types of testing. It requires effective reference temperatures $u_\mathrm{melt}$, $u_\mathrm{surf}$. Furthermore, the weights and the secondary path $\mathcal{C}^\mathrm{wd}$ need to be carefully chosen. Work related to this kind of testing has not been detailed here.

\section*{Acknowledgements}
\noindent This work was supported by the Swedish Foundation for Strategic Research under the contract ``Industrial PhD 2015 -- ID15-0058''.

\bibliographystyle{plain}
\bibliography{bibliography}

\begin{thebibliography}{10}

\bibitem{doi:10.1137/0916069}
Richard~H. Byrd, Peihuang Lu, Jorge Nocedal, and Ciyou Zhu.
\newblock A limited memory algorithm for bound constrained optimization.
\newblock {\em SIAM Journal on Scientific Computing}, 16(5):1190--1208, 1995.

\bibitem{cao}
X.~Cao and B.~Ayalew.
\newblock Partial differential equation-based multivariable control input
  optimization for laser-aided powder deposition processes.
\newblock {\em ASME. J. Manuf. Sci. Eng.}, 138(3):031001--031001--8, 2015.

\bibitem{eagerandtsai}
T.~W. Eagar and N.-S. Tsai.
\newblock Temperature fields produced by traveling distributed heat sources.
\newblock {\em Weld. Res. Suppl.}, 62:346--355, 1983.

\bibitem{forslund}
R.~Forslund, A.~Snis, and S.~Larsson.
\newblock Analytical solution for heat conduction due to a moving {G}aussian
  heat flux with piecewise constant parameters.
\newblock {\em Appl. Math. Model.}, 66, 2019.

\bibitem{garg}
A.~Garg, K.~Tai, and M.M. Savalani.
\newblock State-of-the-art in empirical modelling of rapid prototyping
  processes.
\newblock {\em Rapid Prototyping J.}, 20(2):164--178, 2014.

\bibitem{gong}
H.~Gong, K.~Rafi, T.~Starr, and B.~Stucker.
\newblock The effects of processing parameters on defect regularity in
  {T}i-6{A}l-4{V} parts fabricated by selective laser melting and electron beam
  melting.
\newblock In {\em Solid Freeform Fabrication Symposium}, pages 424--439, 2013.

\bibitem{heinl}
Peter Heinl, Lenka M{\"u}ller, Carolin K{\"o}rner, Robert~F. Singer, and
  Frank~A. M{\"u}ller.
\newblock Cellular {T}i–6{A}l–4{V} structures with interconnected macro
  porosity for bone implants fabricated by selective electron beam melting.
\newblock {\em Acta Biomaterialia}, 4(5):1536 -- 1544, 2008.

\bibitem{HINZE2007657}
Michael Hinze and Stefan Ziegenbalg.
\newblock Optimal control of the free boundary in a two-phase {S}tefan problem.
\newblock {\em J. Comput. Phys.}, 223(2):657--684, 2007.

\bibitem{scipy}
Eric Jones, Travis Oliphant, Pearu Peterson, et~al.
\newblock {SciPy}: Open source scientific tools for {Python}, 2001--.

\bibitem{kamath}
Chandrika Kamath, Bassem El-dasher, Gilbert~F. Gallegos, Wayne~E. King, and
  Aaron Sisto.
\newblock Density of additively-manufactured, 316{L SS} parts using laser
  powder-bed fusion at powers up to 400 {W}.
\newblock {\em Int. J. Adv. Manuf. Technol.}, 74(1):65--78, Sep 2014.

\bibitem{khairallah}
Saad~A. Khairallah, Andrew~T. Anderson, Alexander Rubenchik, and Wayne~E. King.
\newblock Laser powder-bed fusion additive manufacturing: Physics of complex
  melt flow and formation mechanisms of pores, spatter, and denudation zones.
\newblock {\em Acta Materialia}, 108:36--45, 2016.

\bibitem{kingkhairallah}
W.~King, A.~T. Anderson, R.~M. Ferencz, N.~E. Hodge, C.~Kamath, and S.~A.
  Khairallah.
\newblock Overview of modelling and simulation of metal powder bed fusion
  process at {L}awrence {L}ivermore {N}ational {L}aboratory.
\newblock {\em Mater. Sci. Tech. Ser.}, 31(8):957--968, 2015.

\bibitem{kingetal}
W.~E. King, A.~T. Anderson, R.~M. Ferencz, N.~E. Hodge, C.~Kamath, S.~A.
  Khairallah, and A.~M. Rubenchik.
\newblock Laser powder bed fusion additive manufacturing of metals; physics,
  computational, and material challenges.
\newblock {\em Appl. Phys. Rev.}, 2(4), 2015.

\bibitem{Lewis2013}
Adrian~S. Lewis and Michael~L. Overton.
\newblock Nonsmooth optimization via quasi-{N}ewton methods.
\newblock {\em Mathematical Programming}, 141(1):135--163, Oct 2013.

\bibitem{ma}
L.~Ma, J.~Fong, B.~Lane, S~Moylan, J.~Filliben, A.~Heckert, and L.~Levine.
\newblock Using design of experiments in finite element modeling to identify
  critical variables for laser powder bed fusion.
\newblock In {\em Solid Freeform Fabrication Symposium}, pages 219--228, 2015.

\bibitem{john}
J.~Bondestam Malmberg and M.~Wallen{\aa}s.
\newblock Solving the heat equation in connection with electron beam melting.
\newblock Master's thesis, Department of Mathematical Sciences, Mathematics,
  Chalmers University of Technology, 2012.

\bibitem{mani}
M.~Mani, B.~Lane, M.~A. Donmez, S.~Feng, S.~Moylan, and R.~Fesperman.
\newblock Measurement science needs for real-time control of additive
  manufacturing powder bed fusion processes.
\newblock Technical report, National Institute of Standards and Technology,
  Gaithersburg, MD, NIST Interagency/Internal Report (NISTIR), 2015.

\bibitem{marklandkorner}
M.~Markl and C.~K{\"o}rner.
\newblock Multi-scale modeling of powder-bed-based additive manufacturing.
\newblock {\em Annu. Rev. Mater. Res.}, 46:93--123, 2016.

\bibitem{mukherjee}
T.~Mukherjee, J.S. Zuback, A.~De, and T.~DebRoy.
\newblock Printability of alloys for additive manufacturing.
\newblock {\em Sci. Rep.}, 6, 2016.

\bibitem{ning}
Y.~Ning, J.~Y.~H. Fuh, Y.~S. Wong, and H.~T. Loh.
\newblock An intelligent parameter selection system for the direct metal laser
  sintering process.
\newblock {\em Int. J. Prod. Res.}, 42(1):183--199, 2004.

\bibitem{smith}
C.J. Smith, F.~Derguti, E.~Hernandez Nava, M.~Thomas, S.~Tammas-Williams,
  S.~Gulizia, D.~Fraser, and I.~Todd.
\newblock Dimensional accuracy of {E}lectron {B}eam {M}elting ({EBM}) additive
  manufacture with regard to weight optimized truss structures.
\newblock {\em J. Mater. Process Tech.}, 229:128--138, March 2016.

\bibitem{snis2}
A.~Snis.
\newblock Method for production of a three-dimensional body, July~7 2015.
\newblock {US} Patent 9,073,265 B2.

\bibitem{vandenbroucke}
B.~Vandenbroucke and J.~Kruth.
\newblock Selective laser melting of biocompatible metals for rapid
  manufacturing of medical parts.
\newblock {\em Rapid Prototyping J.}, 13(4):196--203, 2007.

\bibitem{Volkov2009}
Oleg Volkov, Bartosz Protas, Wenyuan Liao, and Donn~W. Glander.
\newblock Adjoint-based optimization of thermo-fluid phenomena in welding
  processes.
\newblock {\em J. Eng. Math.}, 65(3):201--220, Nov 2009.

\bibitem{zeng}
K.~Zeng, D.~Pal, and B.~Stucker.
\newblock A review of thermal analysis methods in laser sintering and selective
  laser melting.
\newblock In {\em Solid Freeform Fabrication Symposium}, pages 796--814, 2012.

\bibitem{Zhu:1997:ALF:279232.279236}
Ciyou Zhu, Richard~H. Byrd, Peihuang Lu, and Jorge Nocedal.
\newblock Algorithm 778: {L-BFGS-B}: Fortran subroutines for large-scale
  bound-constrained optimization.
\newblock {\em ACM Trans. Math. Softw.}, 23(4):550--560, December 1997.

\end{thebibliography}

\end{document}